\documentclass[leqno]{amsart}
\usepackage{amsmath,amsfonts,amsthm,amssymb,indentfirst}

\newcommand{\<}{\kern.0833em}

\newtheorem{theorem}{Theorem}
\newtheorem{lemma}[theorem]{Lemma}
\newtheorem{corollary}[theorem]{Corollary}
\newtheorem{proposition}[theorem]{Proposition}
\newtheorem{definition}[theorem]{Definition}
\newtheorem{question}[theorem]{Question}

\newcommand{\V}{\mathbf{V}}
\newcommand{\Pv}{\mathbf{P}}
\newcommand{\X}{\mathbf{X}}
\newcommand{\Su}{\mathbb{S}\,}
\renewcommand{\Pr}{\mathbb{P}\,}
\newcommand{\strt}[1][1.7]{\vrule width0pt height0pt depth#1pt}
\newcommand{\cP}{\raisebox{.12em}{$\<\scriptscriptstyle\coprod$}\kern -.2em\nolinebreak[2]^{\strt}} 
\newcommand{\tto}{\rightarrowtail}

\begin{document}

\begin{center}
\texttt{Comments, corrections, and related references welcomed, as
  always!}\\[.5em]
{\TeX}ed \today
\vspace{2em}
\end{center}

\title{On coproducts in varieties, quasivarieties and prevarieties}%
\thanks{
arXiv reference: 0806.1750\,.
After publication of this note, updates, errata, related references
etc., if found, will be recorded at
http://math.berkeley.edu/%
{$\!\sim$}gbergman\protect\linebreak[0]/papers\,.
}

\keywords{Coproduct of algebras in a variety,
quasivariety or prevariety.
Free algebra on $n$ generators which contains a subalgebra
free on ${>}\,n$ generators.
Amalgamation property.
Number of algebras needed to generate a quasivariety or prevariety.
Symmetric group on an infinite set.}

\author{George M. Bergman}
\address{Department of Mathematics\\
University of California\\
Berkeley, CA 94720-3840\\
USA}
\email{gbergman@math.berkeley.edu}

\begin{abstract}
If the free algebra $F$ on one generator in a variety $\V$
of algebras (in the sense of universal algebra)
has a subalgebra free on two generators, must it also have a
subalgebra free on three generators?
In general, no; but yes if $F$ generates the variety $\V.$

Generalizing the argument, it is shown that if we are given an algebra
and subalgebras, $A_0\supseteq\dots\supseteq\nolinebreak A_n,$
in a prevariety $\!(\Su\Pr\!$-closed class of algebras) $\Pv$ such that
$A_n$ generates $\Pv,$ and also subalgebras
$B_i\subseteq\nolinebreak A_{i-1}$ $(0<i\leq n)$ such that
for each $i>0$ the subalgebra of $A_{i-1}$ generated by $A_i$ and $B_i$
is their coproduct in $\Pv,$ then the subalgebra of $A$ generated by
$B_1,\dots,B_n$ is the coproduct in $\Pv$ of these algebras.

Some further results on coproducts are noted:

If $\Pv$ satisfies the amalgamation property, then one has the stronger
``transitivity'' statement, that if $A$ has a finite family of
subalgebras $(B_i)_{i\in I}$ such that the subalgebra of $A$ generated
by the $B_i$ is their coproduct, and each $B_i$ has a
finite family of subalgebras $(C_{ij})_{j\in J_i}$ with the same
property, then the subalgebra of $A$ generated by all the $C_{ij}$
is their coproduct.

For $\Pv$ a residually small prevariety or an arbitrary
quasivariety, relationships are proved between
the least number of algebras needed to generate $\Pv$ as a prevariety
or quasivariety, and behavior of the coproduct operation in $\Pv.$

It is shown by example that for $B$ a subgroup of the group
$S=\mathrm{Sym}(\Omega)$ of all permutations of an infinite set
$\Omega,$ the group $S$ need not have a subgroup
isomorphic over $B$ to the coproduct with amalgamation $S\,\cP_B\,S.$
But under various additional hypotheses on $B,$
the question remains open.
\end{abstract}
\maketitle

\section{Prologue, for the non-expert.}\label{S.prologue}
It is well known that the free group on two generators contains
a subgroup free on three generators.
Can one deduce, from this alone, that it contains a subgroup free on
four generators?

It is unfair to say ``from this alone'' without
indicating what facts about groups are to be taken for granted.
So suppose we want to use only the fact that groups form a {\em variety}
of algebras in the sense of universal algebra -- a class of
structures consisting of all sets with a family of operations of
specified arities, satisfying a specified list of identities.
Then we can ask, for $\V$ any variety, and $n$ any positive integer,
\begin{equation}\begin{minipage}[c]{35pc}\label{d.freeq}
If in $\V$ the free algebra on $n$ generators
has a subalgebra free on $n+1$ generators, must it have a
subalgebra free on $n+2$ generators?
\end{minipage}\end{equation}
Our first result will be a negative answer to this question,
in the most extreme case, $n=1.$

On the other hand, a fact that is second nature to combinatorial
group theorists is that if $G_1$ and $G_2$ are overgroups of a
common group $H,$ and one forms $G_1\cP_H\,G_2,$
their coproduct with amalgamation of $H$
(in group-theorists' notation and language $G_1 *_H G_2,$ their
free product with amalgamation of $H),$ then the canonical maps of
$G_1$ and $G_2$ into that coproduct are embeddings.
This says that the variety of all groups has
``the amalgamation property''; and
we shall see in \S\ref{S.amalg} that if a variety $\V$ has
this property, then it also has the property that for any algebras
$A_1,\ A_2$ in $\V$ and subalgebras
$B_1\subseteq A_1,\ B_2\subseteq A_1,$ the coproduct
$A_1\cP\,A_2$ contains the coproduct $B_1\cP\,B_2.$
From this it is not hard to show that for any such
$\V,$~(\ref{d.freeq}) has an affirmative answer.

However, the amalgamation property is relatively rare.
For instance, though it is satisfied by the variety of
all groups, and by all varieties of abelian groups, it does
not seem to be satisfied by most other varieties of groups -- in fact,
it is a longstanding open question whether it is satisfied
by any variety of groups other than those
just mentioned \cite[Problem~6]{HN},~\cite[p.422]{HNq}.

But in section~\S\ref{S.free}, after finding our
counterexample to~(\ref{d.freeq}), we shall see
that a different condition,
more common than the amalgamation property, implies a
positive answer to~(\ref{d.freeq});
namely, that the free algebra of rank~$n$ in $\V$
{\em generate} $\V,$ i.e., not lie in any proper subvariety thereof.
(As, for example, the free group of rank~$2$ generates the
variety of all groups.)

In \S\S\ref{S.quasi}-\ref{S.P-indp} we shall generalize this
to a result about when an algebra (not necessarily free) containing
a coproduct of subalgebras, some of which in turn contain coproducts
of subalgebras, will itself contain an ``obvious'' iterated coproduct.
The condition that a
certain one of our algebras generate the class we are working in will
again be a key assumption; not, this time, for generation as
a variety, but as a {\em prevariety,} which means, roughly,
a class of algebras determined by identities
and universal {\em implications}.
(For example, the class of torsion-free groups, i.e., groups
satisfying $(\forall\,x)\ x^n=1\implies x=1$ for each $n>0,$ is
a prevariety.)
The definition of prevariety, and of the related concept of
{\em quasivariety}, are recalled in~\S\ref{S.quasi}.

We end with some further results on
quasivarieties and prevarieties, and a brief
section on subgroups of infinite symmetric groups.

I am indebted to the referee for several helpful suggestions,
and to the editorial staff of the journal for
requesting that I write this introduction for the general reader.

Carrying that suggestion further, I have added,
as \S\ref{S.glossary}, a quick summary of some common terminology
which should make this note readable (if not
light reading) by anyone for whom this prologue was.
Readers not familiar with the basic language of universal algebra
might start with that section.

\section{Free subalgebras of free algebras.}\label{S.free}

The original question that led to this investigation,
\cite[Question~4.5]{embed}, was whether an algebra $A$ in a
variety $\V$ which contains a subalgebra isomorphic to the coproduct
in $\V$ of two copies of itself, $A\cP_\V\,A,$ must also contain a copy
of the $\!3\!$-fold coproduct $A\cP_\V\,A\cP_\V\,A.$
As indicated above, this can fail even for $A$ free of rank~$1:$
a free algebra of rank~$1$ in a variety $\V$ may have a subalgebra
free of rank~$2$ but fail to have any subalgebra free of rank~$3.$
Let us begin by examining how we might concoct such an example.

To do so, we must ``foil'' the obvious ways one would expect
a free $\!3\!$-generator subalgebra to arise.
If $\langle x\rangle$ is free on $x$ and contains a subalgebra
$\langle y,\,z\rangle$ free on $y$ and $z,$ then
$y=px$ and $z=qx$ for some derived unary operations $p,\,q$ of $\V.$
Since $\langle qx\rangle$ is isomorphic to
$\langle x\rangle,$ its subalgebra
corresponding to $\langle y,\,z\rangle,$ namely
$\langle\<pqx,\,qqx\rangle,$ will be free on those two
generators, and one might expect $\langle\<px,\,pqx,\,qqx\rangle$ to be
free on the three indicated generators.
(If it seems to the reader that
it {\em must} be free on those elements, he or she may be
implicitly assuming that $\V$ has the amalgamation property, to
be discussed in~\S\ref{S.amalg}.)

For this to fail, there must be some ternary relation $T$ in the
operations of $\V$ such that $T(px,\,pqx,\,qqx)$ is an identity
in one variable $x,$ but $T(u,\,v,\,w)$ is
not satisfied by all $\!3\!$-tuples of elements of algebras in $\V.$
On the other hand, since $y=px$ and $z=qx$ generate a free
algebra, the relation $T(px,\,pqx,\,qqx)$ implies that
$T(y,pz,qz)$ {\em is} an identity in two variables $y$ and $z$ in $\V.$

Let us pause to note that if we construct a variety with such an
identity $T,$ we will have eliminated one possibility for a free
$\!3\!$-generator subalgebra of $\langle x\rangle$ of rank~$3;$
but every $\!3\!$-tuple of expressions obtained
from $y$ and $z$ using the operations of $\V$ represents
another potential generating set for a free subalgebra.
In principle, we might use different relations to exclude different
$\!3\!$-tuples; but let us see whether we can make do with just one
such relation $T,$ such that
$T(u,v,w)$ holds for all $\!3\!$-tuples $(u,v,w)$ of elements
of $\langle x\rangle.$
Note that in this case, since $\langle x\rangle$ contains a
free algebra of rank two, $\V$ must satisfy
identities saying that $T(u,v,w)$ holds for any
elements $u,\,v,\,w$ of any $\!\V\!$-algebra that
lie in a common $\!2\!$-generator subalgebra.

In testing out this approach, let us temporarily allow structures
involving primitive {\em relations} as well as operations.
Then we could let $\V$ be the class of objects defined by two
primitive unary operations, $p$ and $q,$ and one primitive ternary
relation, $T,$ subject only to the countable family of ``identities''
\begin{equation}\begin{minipage}[c]{35pc}\label{d.T0}
$T(a(y,z),\,b(y,z),\,c(y,z)),$
\end{minipage}\end{equation}
one for each $\!3\!$-tuple of words
$a,\,b,\,c$ in two variables $y,\,z$ and the operations $p,\,q.$
(Of course, since $p$ and $q$ are unary, each of $a,\ b,\ c$ really
just involves one of $y$ or $z.)$
In an object of $\V$ generated by $\leq 2$ elements, $T$ thus holds
identically, so in describing the
structures of $\!\leq 2\!$-generator objects, we can
ignore the relation $T,$ and simply specify the actions of $p$ and $q.$
Since the family of identities~(\ref{d.T0}) by which we have
defined $\V$ includes no identities in the
operations $p$ and $q$ alone, the possible structures of
such objects are simply the structures of $\!M\!$-set, for $M$ the
free monoid on generators $p$ and $q.$
In this monoid $M,$ the left ideal generated by $p$ and $q$ is
free on those two elements; hence in the free $\!\V\!$-object
on one generator $x,$ the elements $px$ and $qx$ satisfy no relations
in $p$ and $q;$ so with $T,$ as noted, also contributing no
information, $px$ and $qx$ indeed generate a free subobject.
On the other hand, if we take the free $\!M\!$-set on three
generators $x,\,y,\,z,$ and define $T$ to hold precisely on those
$\!3\!$-tuples thereof in which all three components lie in a
subalgebra generated by two elements, we see that this satisfies the
definition of a free $\!\V\!$-object on three generators,
and that $T(x,y,z)$ does {\em not} hold.
Hence the free object on one generator does not contain a copy
of the free object on three generators.

Let us now try to mimic the above behavior in a variety of genuine
algebras.
In addition to two unary operations $p$ and $q,$ let us
introduce a $\!0\!$-ary operation $0$ and a ternary operation $t,$
with the idea that the relation
$T(u,v,w)$ will be the condition $t(u,v,w)=0.$
To keep our new operations from complicating our structures more
than necessary, let us introduce some ``nonproliferation'' identities:
\begin{equation}\begin{minipage}[c]{35pc}\label{d.p=q=0}
$p\<0\ =\ q\<0\ =\ p\,t(x,y,z)\ =\ q\,t(x,y,z)\ =\ 0.$
\end{minipage}\vspace{-0.4em}\end{equation}
\begin{equation}\begin{minipage}[c]{35pc}\label{d.t=0}
$t(u,v,w)\ =\ 0$ \ whenever any of $u,v,w$ is either $0,$ or is
itself of the form $t(u',v',w').$
\end{minipage}\end{equation}
Finally, we impose the identities corresponding to~(\ref{d.T0}):
\begin{equation}\begin{minipage}[c]{35pc}\label{d.T}
$t(a(x,y),\<b(x,y),\<c(x,y))\ =\ 0$
\ for all derived operations $a,\ b,\ c$ in two variables.
\end{minipage}\end{equation}

In a free $\!\V\!$-algebra, the elements $t(u,v,w)$ that
are not $0$ may be thought of as ``tags'', showing that certain
$\!3\!$-tuples $(u,v,w)$ obtained from the generators using
$p$ and $q$ alone do not have the form indicated in~(\ref{d.T}).
By~(\ref{d.p=q=0}) and~(\ref{d.t=0}), these elements have essentially
no other effect.
By the same reasoning as
for structures with a primitive relation $T,$ we get

\begin{proposition}\label{P.2_not_3}
Let $\V$ be the variety defined by a zeroary operation $0,$ two
unary operations $p$ and $q,$ and a ternary operation $t,$ subject
to identities \textup{(\ref{d.p=q=0}), (\ref{d.t=0}), (\ref{d.T})}.

Then in the free $\!\V\!$-algebra $F_\V(x)$ on one generator $x,$
the subalgebra generated by $px$ and $qx$ is free on those generators;
but $F_\V(x)$ \textup{(}and hence also the free algebra
on two generators\textup{)} has no subalgebra free on three or
more generators.\qed
\end{proposition}

The above result was based on $F_\V(x)$ satisfying an identity
(namely, $t(x,y,z)=0)$ that did not hold in all of $\V;$
and we might hope that if
$\V$ is a variety where this does not happen, but which, as above, has
unary derived operations $p$ and $q$ such that $px$ and $qx$ are free
generators of the subalgebra $\langle\<px,qx\rangle\subseteq F_\V(x),$
then the subalgebra $\langle\<px,pqx,qqx\rangle$
will have to be free on the indicated three generators.
To investigate this question, consider a ternary relation $T$ in the
operations of $\V$ about which we now merely
assume that $T(px,pqx,qqx)$ holds in $F_\V(x),$
and let us see whether we can deduce that $T$ holds for all
$\!3\!$-tuples of elements of $F_\V(x).$

As noted earlier, the conditions that $T(px,pqx,qqx)$ holds in
$F_\V(x),$ and that $y=px$ and $z=qx$ generate a free algebra
$\langle y,z\rangle,$ show that in that free algebra,
$T(y,pz,qz)$ holds, hence that in any $\!\V\!$-algebra, $T$
holds on any $\!3\!$-tuple whose last two terms are obtained
from a common element by applying $p,$ respectively, $q$ to it.
Let us now apply this observation to a $\!3\!$-tuple in $F_\V(x)$ of
the form $(a(px,qx),px,qx)$ where $a$ is any derived operation of $\V,$
and use the independence of $px$ and $qx$ a second time.
We conclude that $T(a(y,z),y,z)$ holds for every such $a.$
In other words, in any $\!\V\!$-algebra, $T$ holds on every
$\!3\!$-tuple whose first term lies
in the subalgebra generated by the last two terms.

But there is no evident way to carry this process further.
And in fact, we can again get a negative result by the same
technique of realizing $T$ as $t(u,v,w)=0,$
embodying the conditions that we have found $T$ must
satisfy this time, in the system of identities:
\begin{equation}\begin{minipage}[c]{35pc}\label{d.u,pv,qv}
$t(u,\,pv,\,qv)\ =\ 0$ \ for all $u,\,v.$
\end{minipage}\vspace{-0.4em}\end{equation}
\begin{equation}\begin{minipage}[c]{35pc}\label{d.auv,u,v}
$t(a(u,v),\,u,\,v)\ =\ 0$ \ for all $u,\,v,$ and all binary
terms $a.$
\end{minipage}\end{equation}

The one tricky point is to show that the variety so defined now has
the property that there
are no identities satisfied by the free algebra on one generator that
are not identities of the whole variety.
In contrast to the earlier example, our development has not called
on any such identities; but neither has it shown that none exist.
With some work, one can prove this; but an easier approach,
which we will follow, is to let $\V_0$ denote the variety defined by
the identities discussed above, and let our $\V$ be the subvariety
of $\V_0$ generated by the free algebra on one generator therein.
Here are the details.

\begin{proposition}\label{P.not_px,pqx,qqx}
Let $\V_0$ be the variety defined by a zeroary operation $0,$ two unary
operations $p$ and $q,$ a ternary operation $t,$ and the
identities \textup{(\ref{d.p=q=0}), (\ref{d.t=0}), (\ref{d.u,pv,qv}),}
and \textup{(\ref{d.auv,u,v});} and let $\V$ be the subvariety of $\V_0$
generated by the free algebra $F_{\V_0}(x)$ on one generator.
Thus, $F_{\V}(x)=F_{\V_0}(x),$ so $\V$ is generated by $F_{\V}(x).$

In this situation, the subalgebra
$\langle\<px,qx\rangle\subseteq F_{\V}(x)$ is free in $\V$
\textup{(}and in fact in $\V_0)$ on the two generators $px$ and $qx;$
but the subalgebra $\langle\<px,pqx,qqx\rangle$ is not free in
$\V$ \textup{(}and hence not in $\V_0)$ on $px,\ pqx$ and $qqx.$
\end{proposition}

\begin{proof}
The last sentence of the first paragraph is clear in the general
context of a subvariety generated by a free algebra in any variety.

We shall next show that $\langle\<px,qx\rangle\subseteq F_{\V}(x)$
is free on $px$ and $qx$ in $\V_0.$
Since it is a subalgebra of $F_{\V}(x)$ and
hence belongs to $\V,$ it will then {\em a fortiori}
be free on those generators in that subvariety.

To do this, we need to prove that any relation satisfied in $F_{\V}(x)$
by $px$ and $qx$ also holds between $y$ and $z$ in $F_{\V_0}(y,z).$
Let $M$ again denote the free monoid on the two symbols $p$ and $q.$
We know as before that the elements of $F_{\V}(x)$ obtained
from $x$ using $p$ and $q$ alone form a free
$\!M\!$-set on one generator, and hence that the sub-$\!M\!$-set
$M\{px,qx\}$ is free as an $\!M\!$-set on $px$ and $qx.$
Thus, if the given relation satisfied
by $px$ and $qx$ involves only the operations
$p$ and $q,$ it will indeed be satisfied by $y$ and $z.$
Hence in what follows, we may assume the relation involves
$t$ and/or $0.$

Now by~(\ref{d.p=q=0}) and~(\ref{d.t=0}), if either side of our
relation in $px$ and $qx$
involves $0$ or $t$ other than in the outermost position,
that side equals $0,$ and the corresponding expression in $y$ and
$z$ does as well; so we can replace that side by $0$ in our relation.
Given the forms of the identities~(\ref{d.u,pv,qv})
and~(\ref{d.auv,u,v}), it is not hard to see that to complete
our proof it will suffice to show that if
\begin{equation}\begin{minipage}[c]{35pc}\label{d.twww}
$t(b(px,qx),\,c(px,qx),\,d(px,qx))\ =\ 0$
\end{minipage}\end{equation}
is an identity of $\V_0,$ where
\begin{equation}\begin{minipage}[c]{35pc}\label{d.bcd}
$b(px,qx),\ c(px,qx),\ d(px,qx)\,\in\,M\{px,qx\},$
\end{minipage}\end{equation}
then $\V_0$ also satisfies the identity
\begin{equation}\begin{minipage}[c]{35pc}\label{s.tbcd}
$t(b(y,z),\,c(y,z),\,d(y,z))=0.$
\end{minipage}\end{equation}
Moreover,~(\ref{d.p=q=0}) and~(\ref{d.t=0}) yield no
relations of the form~(\ref{d.twww}) satisfying~(\ref{d.bcd}),
so we need only look at relations~(\ref{d.twww}) of
the forms~(\ref{d.u,pv,qv}) and~(\ref{d.auv,u,v}).

An instance of~(\ref{d.auv,u,v}) can have the form~(\ref{d.twww})
only if the given $u$ and $v$ have the forms $c(px,qx)$ and $d(px,qx);$
but then putting $y$ and $z$ in place of $px$ and $qx$ in that instance
of~(\ref{d.auv,u,v}) again gives an instance of~(\ref{d.auv,u,v}),
and hence a relation in $F_{\V_0}(y,z),$ as required.
If an instance of~(\ref{d.u,pv,qv}) has the
form~(\ref{d.twww}), then we have $u=b(px,qx),$ but
there are two possibilities for the element $v:$
it can either be $x,$ or of the form $e(px,qx).$
In the former case, this instance of~(\ref{d.u,pv,qv}) is also an
instance of~(\ref{d.auv,u,v}), and the preceding argument applies.
In the latter case, the relation has the
form $t(b(px,qx),\,pe(px,qx),\,qe(px,qx))=0,$ and we see
that $t(b(y,z),\,pe(y,z),\,qe(y,z))=0$ is again an instance
of~(\ref{d.u,pv,qv}), and hence a relation in $F_{\V_0}(y,z).$
This completes the proof that $\langle\<px,qx\rangle$ is free
on $px$ and $qx.$

To see, finally, that $\langle\<px,pqx,qqx\rangle$
is not free on the indicated generators in $\V,$
we note that $F_{\V_0}(x),$ which generates $\V,$ has $\!3\!$-tuples
of elements of $M\{x\}$
to which neither~(\ref{d.u,pv,qv}) nor~(\ref{d.auv,u,v}) applies,
e.g., $(x,qx,px).$
Hence $t(x,y,z)=0$ is not an identity of $\V;$
hence the elements $px,$ $pqx,$ $qqx,$ which {\em do} satisfy
that relation, cannot be free generators of a free subalgebra.
\end{proof}

After obtaining the above result, I wondered whether for {\em every}
$\!3\!$-tuple $(ax,bx,cx)$ in $M\{x\},$ one could find a ternary
relation $T_{a,b,c}$ on $M\{x\}$ that could be embodied in a
construction like the above, giving an algebra in which
that $\!3\!$-tuple was not a free generating set.
If so, then it would seem that by defining a variety with
operations $0,\ p,$ and $q$ and countably
many ternary operations $t_{a,b,c},$ one for
each such choice of $a,$ $b$ and $c,$ one should be able to get an
example where, as above, $F_{\V}(x)$ generated $\V$ and
$\langle\<px,qx\rangle\subseteq F_\V(x)$ was free on
$px,\ qx,$ but where $F_\V(x)$ contained no subalgebra free on three
generators.

But just a bit more experimentation revealed $\!3\!$-tuples $(ax,bx,cx)$
for which no $T_{a,b,c}$ with the desired property exists.
Translating the resulting obstruction into a proof
of a positive statement, this is

\begin{proposition}\label{P.2->3}
Let $\V$ be a variety of algebras such that the free
algebra $F_\V(x)$ on one generator generates $\V$ as a variety,
and contains a subalgebra free of rank $2$ in $\V,$ say
on generators $px$ and $qx,$ where
$p$ and $q$ are derived unary operations of $\V.$
Then the subalgebra
$\langle\<px,\,pqx,\,pqqx\rangle\subseteq F_\V(x)$
is free in $\V$ on the indicated three generators.
\end{proposition}

\begin{proof}
It will suffice to show that for any three elements
\begin{equation}\begin{minipage}[c]{35pc}\label{d.abc}
$ax,\ bx,\ cx\in F_\V(x),$
\end{minipage}\end{equation}
there exists a homomorphism
$\langle\<px, pqx,pqqx\rangle\rightarrow F_\V(x)$ carrying
$px,\ pqx,\ pqqx$ to $ax,\ bx,\ cx$ respectively, since this
will show that every relation satisfied by $px,\ pqx$ and $pqqx$
is an identity of $F_\V(x),$ and hence, by hypothesis, of~$\V.$

Given elements~(\ref{d.abc}), let us first use the freeness of
$\langle\<px,qx\rangle$ to get a homomorphism
$f:\langle\<px,qx\rangle\to F_\V(x)$ carrying
$px$ to $aqqx,$ and $qx$ to $x.$
Thus the image of $(px,\ pqx,\ pqqx)$ under this map
is $(aqqx,\ px,\ pqx).$
Since this $\!3\!$-tuple, and hence the subalgebra it
generates, again lies in $\langle\<px,qx\rangle,$ we
can compose this homomorphism with another homomorphism,
$g:\langle\<px,qx\rangle\to F_\V(x);$
let this take $px$ to $bqx$ and $qx$ to $x.$
This takes the preceding $\!3\!$-tuple to $(aqx, bqx, px).$
Finally, mapping $\langle\<px,qx\rangle$ to $F_\V(x)$
by the homomorphism $h$ sending $px$ to $cx$ and $qx$ to $x,$ we get
the desired $\!3\!$-tuple $(ax, bx, cx).$
Hence, the composite $h\<gf:\langle\<px,\,pqx,\,pqqx\rangle\to F_\V(x)$
acts as required.
\end{proof}

In the question we have answered, the choice of ranks one, two
and three was, of course, made to
give a concrete test problem.
This, and the restriction to free algebras
rather than coproducts of general algebras,
make our counterexamples, Propositions~\ref{P.2_not_3}
and~\ref{P.not_px,pqx,qqx}, formally stronger, but our positive result,
Proposition~\ref{P.2->3}, weaker than the corresponding result
without those restrictions.
We shall generalize Proposition~\ref{P.2->3} in the next two
sections so as to remove these restrictions.

\section{Prevarieties and quasivarieties.}\label{S.quasi}

In the proof of Proposition~\ref{P.2->3}, we used
the fact that if $\V$ is the variety generated by an
algebra $A,$ then a $\!\V\!$-algebra
generated by a family of elements, $B=\langle\{x_i\mid i\in I\}\rangle,$
is free on those generators if and only if there exist homomorphisms
$B\to A$ taking the $x_i$ to all choices
of $\!I\!$-tuples of elements of $A.$
For our generalization, we would like to say similarly
that if an algebra $B$ is generated by a family of subalgebras
$B_i$ $(i\in I),$ then it is their coproduct if and only if
every system of homomorphisms from the algebras $B_i$ to our
given algebra $A$ extends to a homomorphism $B\to A.$
We shall see that this is true for coproducts, not in the
{\em variety} generated by $A,$ but in the {\em prevariety}
so generated (definition recalled below).

There are a few points of notation and terminology in which
usage is not uniform; we begin by addressing these.

First, we admit the empty algebra when the operations
of our algebras include no zeroary operations.

Second, note that the operators $\mathbb{H}\,,$ $\Su$ and $\Pr$ on
classes of algebras that appear in Birkhoff's Theorem and related
results each come in two slightly different flavors.
One may associate to a class $\X$ of algebras the class of all
{\em factor algebras} of members of $\X$ by congruences, or the class
of algebras {\em isomorphic} to such factor algebras, i.e., the
homomorphic images of members of $\X.$
Likewise, one may associate to $\X$ the family of {\em subalgebras}
of members of $\X,$ or the family of algebras {\em isomorphic} to such
subalgebras; i.e., algebras {\em embeddable} in members of $\X.$
And finally, we may associate to $\X$ the class of {\em direct product}
algebras constructed from members of $\X$ in the standard way as
algebras of tuples, or the class of
algebras {\em isomorphic} to algebras so constructed, i.e., algebras
$P$ that admit a family of maps to the indicated members of $\X$ giving
$P$ the universal property of their direct product.
It is probably an accident of history that the symbols
$\mathbb{H}\,,$ $\Su$
and $\Pr$ were assigned, in two cases (subalgebras and products)
to particular explicit constructions, but in the remaining case
(homomorphic images) to the isomorphism-closed concept.
The standard remedy is to introduce an operator $\mathbb{I}\,,$
taking every class $\X$ of algebras to the class of algebras
isomorphic to members of $\X,$ and apply $\mathbb{I}$ in conjunction
with $\Su$ and $\Pr$ when the wider construction is desired.
But that wider construction usually {\em is} what is desired, so
(following \cite{M+M+T}) we will use the less standard definitions:

\begin{definition}\label{D.HSP}
If $\X$ is a class of algebras of the same type, then
$\Su\X$ will denote the class of algebras isomorphic to
subalgebras of algebras in $\X,$ and
$\Pr\X$ the class of algebras isomorphic to
direct products of algebras in $\X$
\textup{(}including the direct product of the empty
family, the $\!1\!$-element algebra\textup{)}.
As usual, $\mathbb{H}\,\X$ will denote the class of homomorphic
images of algebras in $\X.$
\end{definition}

A third point on which terminology is divided concerns the
definition of ``quasivariety''.
Both usages agree that this means a class of algebras
$A$ determined by a set of conditions of the form
\begin{equation}\begin{minipage}[c]{35pc}\label{d.Horn}
$(\forall~x\in A^I)\ \ (\bigwedge_{j\in J}\ a_j(x)=b_j(x))\ \implies
\ c(x)=d(x),$
\end{minipage}\end{equation}
where $I$ and $J$ are sets, and the $a_j$ and $b_j$ and $c$
and $d$ are $\!I\!$-ary terms in the algebra operations.
$(J$ may be empty, in which case~(\ref{d.Horn}) represents an
ordinary identity.)
The point of disagreement is whether $I$ and $J$ are required to be
finite.
The more standard usage, which, somewhat reluctantly, I will follow,
assumes this; a class of algebras defined by sentences~(\ref{d.Horn})
where $I$ and $J$ are not required to be finite
is then called a {\em prevariety}.
The other usage is exemplified by \cite{JA+LS}, where
``quasivariety'' is defined with no finiteness restriction on
$I$ and $J,$ while ``prevariety'' is used for
a still more general sort of class of algebras (typified by monoids
in which every element is invertible; i.e.\ groups regarded as monoids).

(My discomfort with the standard usage is
that the prefix \mbox{``pre-''} suggests a
concept used mainly for technical purposes in the development of
another concept, as in ``preorder'', ``presheaf'' and ``precategory''.
Also, the relationship between ``prevariety'' and ``quasivariety''
is not
mnemonic, as ``quasivariety'' and ``elementary quasivariety'' would be.
Incidentally, if $\X$ is a finite set of finite algebras,
the prevariety and the quasivariety that it generates are the same, so
works like \cite{DC+BD} don't have to distinguish the concepts.)

The concept of quasivariety is a natural one only for finitary
algebras.
(The constructions of reduced products and ultraproducts, occurring
in standard characterizations of quasivarieties, are not in general
defined on infinitary algebras.)
Most of our results on prevarieties will not require finitariness;
so algebras comprising prevarieties will not be assumed
finitary unless this is explicitly stated.

We summarize this and some related conventions in

\begin{definition}\label{D.qv}
A {\em prevariety} will mean a class $\Pv$ of algebras of a given
\textup{(}not necessarily finitary\textup{)}
type that can be defined by a class of conditions of the
form\textup{~(\ref{d.Horn})}; equivalently that is
closed under the operators $\Su$ and~$\Pr.$

A prevariety $\Pv$ which is finitary \textup{(}i..e, every
primitive operation of which has finite arity\textup{),}
will be called a {\em quasivariety} if
\begin{equation}\begin{minipage}[c]{35pc}\label{d.Dqv_Horn}
$\Pv$ can be defined by conditions~\textup{(\ref{d.Horn})}
in each of which $I$ and $J$ are finite,
\end{minipage}\end{equation}
equivalently,
\begin{equation}\begin{minipage}[c]{35pc}\label{d.Dqv_ult}
$\Pv$ is closed under ultraproducts,
\end{minipage}\end{equation}
equivalently,
\begin{equation}\begin{minipage}[c]{35pc}\label{d.Dqv_red}
$\Pv$ is closed under reduced products.
\end{minipage}\end{equation}

If $\X$ is a class of algebras of a given type, the least
prevariety containing $\X,$ namely, $\Su\Pr\X,$ will be called
the {\em prevariety generated by} $\X.$
Likewise, if the type is finitary, the least quasivariety containing
$\X,$ namely, $\Su\Pr\Pr_\mathrm{\!ult}\,\X=\Su\Pr_\mathrm{\!red}\ \X,$
where $\Pr_\mathrm{\!ult}$ and $\Pr_\mathrm{\!red}$ denote,
respectively, the constructions of ultraproducts and reduced products
\textup{(}and algebras isomorphic thereto\textup{)},
will be called the {\em quasivariety generated by} $\X.$
Again without the assumption of finitariness,
the least variety containing $\X,$ namely, $\mathbb{H}\,\Su\Pr\X,$
will be called the {\em variety generated by} $\X.$
\end{definition}

In any prevariety, one has algebras presented by arbitrary systems of
generators and relations.
In particular, every family of algebras has a coproduct.
A useful characterization of these is

\begin{lemma}\label{L.cPiff}
Let $\X$ be a class of algebras of a given type, let $\Pv=\Su\Pr\X$ be
the prevariety generated by $\X,$ let $B$ be an algebra in $\Pv,$ and
let $f_i: B_i\to B$ $(i\in I)$ be a family
of maps from algebras in $\Pv$ into $B.$

Then the algebra $B$ is a coproduct of the $B_i$ in $\Pv,$ with
the $f_i$ as the coprojection maps, if and only if the following
two conditions are satisfied:
\begin{equation}\begin{minipage}[c]{35pc}\label{d.BgenbyBi}
$B$ is generated as an algebra by the union of the images $f_i(B_i).$
\end{minipage}\end{equation}
\begin{equation}\begin{minipage}[c]{35pc}\label{d.extdtoB}
For every $A$ in our generating class $\X,$ and every choice of
a family of maps $g_i: B_i\to A$ $(i\in I),$ there exists
a homomorphism $g: B\to A$ such that $g_i=g f_i$ for all $i\in I.$
\end{minipage}\end{equation}
\end{lemma}

\noindent
{\it Sketch of proof.}
``Only if'' is straightforward:
the necessity of~(\ref{d.BgenbyBi})
is shown, as usual, by applying the universal property of $B$ as a
coproduct to the maps $f_i,$ regarded as taking the $B_i$ into the
subalgebra $C$ of $B$ that they together generate (which
belongs to $\Pv,$ since $\Pv$ is closed under taking subalgebras);
while the necessity of~(\ref{d.extdtoB})
is a case of that universal property.

Conversely, assuming~(\ref{d.BgenbyBi}) and~(\ref{d.extdtoB}),
let us show that $B$ and the $f_i$ satisfy the universal property
of the coproduct.
Let $C$ be any algebra in $\Pv,$ given with
homomorphisms $a_i: B_i\to C.$

If there exists a homomorphism $a: B\to C$ with $a_i=a f_i$ for
all $i,$ then by~(\ref{d.BgenbyBi}) it will be unique.

To see that such a map exists, we write $C$ as a subalgebra
of a direct product $\prod_{j\in J} A_j$ with all $A_j$ in $\X.$
Then for each $j\in J,$ the composites of the given maps $a_i: B_i\to C$
with the $\!j\!$-th projection $C\to A_j$ give a system of
maps $a_{ij}: B_i\to A_j$ $(i\in I).$
By~(\ref{d.extdtoB}), for each $j$ the $a_{ij}$ are induced by a single
map $a_{*j}:B\to A_j;$ doing this for all $j\in J$ gives a map
$a:B\to\prod_J A_j,$ whose restriction to each $f_i(B_i)\subseteq B$
lies in $C\subseteq\prod_{j\in J} A_j.$
Hence $a(B)$ lies in $C$ by~(\ref{d.BgenbyBi}).
The relations
$a_{ij}=a_{*j} f_i$ now show that $a_i=a f_i,$ as required.
\qed\vspace{.5em}

Remarks: We shall see in \S\ref{S.compat1} that it can
happen that though each $B_i$ lies in $\Pv=\Su\Pr\X,$
no $A\in\X$ simultaneously admits maps from all $B_i.$
In that case, condition~(\ref{d.extdtoB}) is vacuous, and the lemma says
that~(\ref{d.BgenbyBi}) characterizes the coproduct $\coprod_\Pv B_i.$
Though
implausible-sounding, this is correct: in that case
an algebra $B$ with maps of the $B_i$ into
it, to belong to $\Pv=\Su\Pr\X,$ must embed in the
product of the vacuous family of members of $\X,$ hence
can have at most one element, so there is hardly any way it
can differ from the desired universal object;~(\ref{d.BgenbyBi}) merely
guarantees that if all $B_i$ are empty, $B$ is also.

In a different direction, taking $I=\emptyset$ in the above
result and recalling that a coproduct of the empty family of
objects in a category is an initial object of the category,
the result says that an algebra $B$ is initial in $\Pv$
if and only if it is generated by the empty set and admits
a homomorphism into each $A\in\X.$

\section{$\!\Pv\!$-independent subalgebras.}\label{S.P-indp}

\begin{definition}\label{D.P-indp}
If $A$ is an algebra in a prevariety $\Pv,$ we shall call
a family of subalgebras $B_i\subseteq A$ $(i\in I)$
{\em $\!\Pv\!$-independent} if the subalgebra $B\subseteq A$ that they
generate, given with the system of inclusion maps $B_i\to B,$ is
a coproduct of the $B_i$ in $\Pv.$
\end{definition}

Here, finally,
is the promised generalization of Proposition~\ref{P.2->3}.

\begin{theorem}\label{T.P-indp}
Suppose that $\Pv$ is a prevariety of algebras and $A_0$ an algebra
in $\Pv,$ and that for some natural number $n$ we are given
subalgebras $A_1,\dots,A_n$ and $B_1,\dots,B_n$ of $A_0,$ such that
for $i=1\dots n,$ $A_i$ and $B_i$ are $\!\Pv\!$-independent,
and are both contained in $A_{i-1}.$
Assume, moreover, that $\Pv$ is generated
as a prevariety
by $A_n.$

Then $B_1,\dots,B_n$ are $\!\Pv\!$-independent.
\end{theorem}

\begin{proof}
Let us prove by induction on $i=0,\dots,n$ a statement a little
stronger than what we will need for $i=n,$ namely that for every
system of homomorphisms $f_j: B_j\to A_i$ $(j=1,\dots,i),$
there exists a unique homomorphism $f$ from the subalgebra of $A_0$
generated by $B_1,\dots,B_i$ and $A_i$ into $A_i$ which acts on each
$B_j$ $(j=1,\dots,i)$ as $f_j,$ and which acts as the identity on $A_i.$

This is clear for $i=0.$
Let $0<i\leq n,$ inductively assume the result
for $i-1,$ and suppose we are given $f_j: B_j\to A_i$ $(j=1,\dots,i).$
Since $A_i\subseteq A_{i-1},$ our inductive hypothesis gives
us a homomorphism $g$ from the subalgebra of $A_0$ generated by
$B_1,\dots,B_{i-1}$ and $A_{i-1}$ into $A_{i-1}$ which agrees
with $f_j$ for $j=1,\dots,i-1,$ and is the identity on $A_{i-1}.$
Note that $g$ will carry the subalgebra generated by
$B_1,\dots,B_i$ and $A_i$ into the subalgebra generated
by $A_i$ (into which it carries $B_1,\dots,B_{i-1}$ and $A_i)$
and $B_i$ (which is contained in $A_{i-1},$ and so
is left fixed).

But by assumption, that subalgebra is the coproduct of
$A_i$ and $B_i,$ so we can map it into $A_i$ by a homomorphism $h$
which acts as the identity on $A_i$ and as $f_i$ on $B_i.$
Now $f=hg$ clearly has the property required for our inductive step.

Taking the $i=n$ case of our result, and ignoring the condition
that $f$ be the identity on $A_n,$ we see that
the subalgebra $B\subseteq A_0$ generated by $B_1,\dots,B_n$
satisfies~(\ref{d.extdtoB}) for $\X$ the singleton family $\{A_n\}.$
Since by assumption $A_n$ generates $\Pv,$ Lemma~\ref{L.cPiff} tells
us that $B$ is the coproduct of the $B_i$ in $\Pv.$
\end{proof}

Remark: We might call a family of subalgebras $B_i$ of an algebra
$A$ in a prevariety $\Pv,$ given with a distinguished member $B_0$
which generates $\Pv,$
``almost $\!\Pv\!$-independent'' if every family of homomorphisms
$f_i: B_i\to B_0$ such that $f_0$ is the identity map of $B_0$ can be
realized by a homomorphism on the subalgebra generated by the $B_i.$
We see from the proof of Theorem~\ref{T.P-indp} that that theorem
remains true if the $\!\Pv\!$-independence
hypothesis is weakened to say that each
pair $(A_i,\,B_i),$ with $A_i$ taken as the distinguished member,
is almost $\!\Pv\!$-independent, and the conclusion
strengthened to say that the $\!n{+}1\!$-tuple
$(A_n,B_1,\dots,B_n),$ with $A_n$ as distinguished member,
is almost $\!\Pv\!$-independent.
The condition of almost $\!\Pv\!$-independence seemed too technical
to use in the formal statement of the theorem; but one might keep it
in mind.
It is interesting that while Proposition~\ref{P.not_px,pqx,qqx}
showed that in the situation described there, the subalgebras
$\langle\<px\rangle,$ $\langle\<pqx\rangle$ and $\langle qqx\rangle$
of $F_\V(x)$ were not $\!\V\!$-independent, the above proof
shows that, with the last of them taken as distinguished,
they are {\em almost} $\!\V\!$-independent.

Note that Theorem~\ref{T.P-indp} holds
even in the case $n=0:$  If $A_0$ generates $\Pv,$ then
the subalgebra of $A_0$ generated by the empty set is
the initial object of $\Pv.$

A case of Theorem~\ref{T.P-indp} with a simpler hypothesis is

\begin{corollary}\label{C.A_i=A}
Suppose $A$ and $B_1,...,B_n$ are algebras in a
prevariety $\Pv,$ such that $A$ generates $\Pv,$ and such that for
each $i,$ $A$ contains an isomorphic copy of
$A\cP_\Pv\,B_i.$
Then $A$ contains an isomorphic copy of~$\coprod_\Pv^{i=1,\dots,n}B_i.$
\qed\vspace{.5em}
\end{corollary}

Recall next that a free algebra in a prevariety $\Pv$ is also free
on the same generators in the variety $\V$ generated by $\Pv.$
Hence we can apply the above results to free algebras in a variety,
and obtain the following result extending Proposition~\ref{P.2->3}
(though we omit, for
brevity, the explicit description of the free generators).

\begin{corollary}\label{C.m<n<N}
Suppose $\V$ is a variety and $m<n$ are positive integers such
that the free $\!\V\!$-algebra $F_\V(x_1,\dots,x_m)$
on $m$ generators has a subalgebra free on $n$ generators,
and such that $\V$ is generated as a variety by $F_\V(x_1,\dots,x_m).$
Then for every natural number $N,$
$F_\V(x_1,\dots,x_m)$ has a subalgebra free on $N$ generators.
\end{corollary}

\begin{proof}
A free $\!\V\!$-algebra $F_\V(x_1,\dots,x_n)$
has subalgebras free on all smaller numbers of generators;
so the above hypothesis implies that $F_\V(x_1,\dots,x_m)$ has
a subalgebra free on $m+1$ generators.
This is a coproduct of a free algebra on $m$ generators and a
free algebra on one generator,
so we get the hypothesis of Corollary~\ref{C.A_i=A} with
$\Pv$ the prevariety generated by $A=F_\V(x_1,\dots,x_m),$
the $n$ of that corollary taken to be $N,$ and
each $B_i$ taken to be free on one generator.
The conclusion shows that $F_\V(x_1,\dots,x_m)$
has a subalgebra free on $N$ generators in the prevariety it generates.
As noted, a free algebra in a
prevariety $\Pv$ is also free in the variety $\V$ generated by $\Pv.$
\end{proof}

Can we strengthen this result to give free subalgebras of
countably infinite rank?
Yes if our algebras are finitary.
We will need

\begin{lemma}\label{L.dirsys}
Let $\Pv$ be a prevariety of {\em finitary} algebras, $A$ an algebra
in $\Pv,$ and $(B_i)_{i\in I}$ a family of subalgebras of $A,$
such that every finite subset $I_0\subseteq I$ is contained
in a subset $I_1\subseteq I$ such that the family of subalgebras
$(B_i)_{i\in I_1}$ is $\!\Pv\!$-independent.
Then $(B_i)_{i\in I}$ is $\!\Pv\!$-independent.
\end{lemma}

\begin{proof}
We need to show that the map $f_I:\coprod_\Pv^I B_i\to A$
whose composite with each coprojection $q_j:B_j\to\coprod_\Pv^I B_i$
is the inclusion of $B_j$ in $A$ is one-to-one.
By finitariness of $\Pv,$ every element of $\coprod_\Pv^I B_i$
lies in the subalgebra generated by finitely many of the $B_i,$
hence it will suffice to show that for any finite
subset $I_0\subseteq I,$ the restriction of $f_I$ to the subalgebra of
$\coprod_\Pv^I B_i$ generated by $\{B_i\mid i\in I_0\}$ is one-to-one.
By assumption, $I_0$ is contained in a subset $I_1$ such that
the family $(B_i)_{i\in I_1}$ is $\!\Pv\!$-independent;
hence the canonical map $f_{I_1}:\coprod_\Pv^{I_1} B_i\to A$
is one-to-one; but that map factors through $f_I,$ so
$f_I$ is one-to-one on its image, the subalgebra of
$\coprod_\Pv^I B_i$ generated by $\{B_i\mid i\in I_1\},$ hence on the
smaller subalgebra generated by $\{B_i\mid i\in I_0\},$ as required.
\end{proof}

(We shall see in \S\ref{S.amalg} and \S\ref{S.compat2}
respectively that if a prevariety $\Pv$ either satisfies
the amalgamation property -- which is not in general the case in the
situation we are interested in here -- or is
generated as a prevariety by a single algebra -- which is
true in the situation to which we are about to apply the above lemma --
then any subfamily of
a $\!\Pv\!$-independent family of subalgebras is $\!\Pv\!$-independent;
so in such cases, the hypothesis of the above lemma can be simplified
merely to say that every finite subset of $I$ is $\!\Pv\!$-independent.
But in a general prevariety $\Pv,$ a subfamily of
a $\!\Pv\!$-independent family need not be $\!\Pv\!$-independent, hence
that simplified statement does not carry the full force of the lemma.
For an example of $\!\Pv\!$-independence not carrying over
to subfamilies, take for $\Pv$ the variety $\V$ of monoids with two
distinguished elements $x$ and $y,$ let $A$ be the $\!\V\!$-algebra
generated by a universal
$\!2\!$-sided inverse to $x,$ denoted $x^{-1},$ let $B_1$ and $B_2$
both be the subalgebra of $A$ generated by $u=x^{-1}y,$ which is a
free monoid on two generators $u$ and $x,$ regarded as a member
of $\V$ by setting $y=xu,$ and let $B_3$ be the whole algebra $A.$
It is not hard to verify that in $B_1\cP_\V\,B_2,$ the images $u_1,$
$u_2$ of the copies of $u$
from $B_1$ and $B_2$ are distinct (though they satisfy $xu_1=xu_2).$
Since in $A$ itself, in contrast, their images are equal, $B_1$ and
$B_2$ are not $\!\V\!$-independent subalgebras of $A.$
But in $B_1\cP_\V\,B_2\cP_\V\,B_3,$ the properties of $\!2\!$-sided
inverses force the generators of $B_1$ and $B_2$ to fall together with
the corresponding elements of $B_3,$ so the family consisting of
these three subalgebras satisfies the definition of
$\!\V\!$-independence.) \vspace{.2em}

Combining the above lemma with our earlier results, we get

\begin{corollary}\label{C.ctbl}
Let $\Pv,$ in~\textup{(i)} and~\textup{(ii)} below, be a prevariety of
{\em finitary} algebras, and $\V,$ in~\textup{(iii)},
a variety of such algebras.
Then\\[.2em]
\textup{(i)} \ If
$A_0\supseteq A_1\supseteq\dots\supseteq A_i\supseteq\dots$
are algebras in $\Pv$ such that every $A_i$ generates $\Pv$ as
a prevariety; and if for each $i>0,$ $B_i$ is a subalgebra of
$A_{i-1}$ such that $A_i$ and $B_i$ are $\!\Pv\!$-independent,
then the countable family $(B_i)_{i>0}$ is $\!\Pv\!$-independent.

Hence,\\[.2em]
\textup{(ii)} \ If $A$ is an algebra which generates $\Pv$
as a prevariety, and we are given a countable family
of algebras $(B_i)_{i>0}$ in $\Pv,$ such that
for each $i,$ $A$ has a subalgebra isomorphic to $A\cP_\Pv\,B_i,$
then $A$ has a subalgebra isomorphic to $\coprod_\Pv^{i>0} B_i.$

Hence,\\[.2em]
\textup{(iii)} \ If for some positive integer $m$ the free algebra
$F_\V(x_1,\dots,x_m)$ generates $\V$ as a variety, and contains
a subalgebra free on $>m$ generators, then it
contains a subalgebra free on countably many generators.
\qed\end{corollary}

Lemma~\ref{L.dirsys} and Corollary~\ref{C.ctbl} both
fail if the assumption that our algebras are finitary is deleted.
To see this for the lemma, let $\V$ be the variety determined by
one operation $a$ of countably infinite arity, and identities saying
that whenever two of $x_0,x_1,\dots$ are equal, we have
\begin{equation}\begin{minipage}[c]{35pc}\label{d.a=x_0}
$a(x_0,\,x_1,\,\dots)\ =\ x_0.$
\end{minipage}\end{equation}
Let $A$ be a countably infinite set, on which $a$ is
defined by letting~(\ref{d.a=x_0}) hold for {\em all} $x_0,x_1,\dots.$
Then every finite subset of $A$ is a free subalgebra on
that set, from which one sees that any finite family of distinct
singleton subsets is an independent set of subalgebras;
but the set of all of these is not independent, because their coproduct
in the variety $\V,$ the free $\!\V\!$-algebra on countably
many generators, does not satisfy~(\ref{d.a=x_0}) identically.

To show that the statements of Corollary~\ref{C.ctbl} all need
the finitariness condition, it suffices to give a counterexample
to statement~(iii) in the absence of that condition.
The idea will be the same as above, but the details
are more complicated, and I will be a little sketchy.

The variety in question will have countably many zeroary
operations $c_0,\,c_1,\dots,$ two unary operations $p$ and $q,$
an operation $a$ of countable arity, and an additional zeroary operation
$0,$ satisfying the analogs of~(\ref{d.p=q=0}) and~(\ref{d.t=0})
with $a$ in place of $t.$
As in \S\ref{S.free}, let $M$ denote the free monoid on the symbols
$p$ and $q.$
Let $\V_0$ be defined by the abovementioned
analogs of~(\ref{d.p=q=0}) and~(\ref{d.t=0}),
together with the (uncountable) family of identities saying that
\begin{equation}\begin{minipage}[c]{35pc}\label{d.a=0}
$a(x_0,x_1,\dots)=0$ if infinitely many of the $x_i$
belong to $M\{u\}$ for some common element $u.$
\end{minipage}\end{equation}

These identities do not imply $a(c_0,c_1,\dots)=0,$ so
$a(x_0,x_1,\dots)=0$ is not an identity in any free algebra in $\V_0.$
Once again, let $\V$ be the subvariety of $\V_0$ generated
by the free algebra $F_{\V}(x)$ on one generator.

One finds that the subalgebra
$\langle px,\,qx\rangle\subseteq F_{\V}(x)=F_{\V_0}(x)$ is free on
$px$ and $qx.$
The key point is that if an element $a(x_0,x_1,\dots)$
with $x_0,x_1,\dots\in\langle px,qx\rangle$ equals $0$ in $F_{\V}(x)$
by an application of~(\ref{d.a=0}),
and the element $u$ of the
hypothesis of~(\ref{d.a=0}) is $x,$ then the infinite family
of elements of $M\{u\}$ in question will be the union of a family of
elements of $M\{px\}$ and a family of elements of $M\{qx\},$
one of which must still be infinite; so the relation
$a(x_0,x_1,\dots)=0$ still follows from the expressions for the
$x_i$ in terms of $px$ and~$qx.$

However, I claim that $F_{\V}(x)$ contains
no subalgebra free on countably many generators.
For note that a family of independent elements of $F_{\V}(x)$ cannot
include the value of any primitive or derived zeroary operation
(since their behavior under homomorphisms is not free), nor
any element obtained with the help of $a,$ by the analogs
of~(\ref{d.p=q=0}) and~(\ref{d.t=0});
hence such a family must lie entirely in $M\{x\}.$
But by~(\ref{d.a=0}) (with
$u=x),$ any infinite family $x_0,\ x_1,\ ,\dots$
of elements of $M\{x\}$ satisfies the relation $a(x_0,x_1,\dots)=0,$
which we have seen is not an identity of $\V;$
so no infinite family of elements of $F_{\V}(x)$ is independent.

\section{Some questions.}\label{S.P-indp-Qs}

Proposition~\ref{P.2_not_3} shows that Corollary~\ref{C.m<n<N}
becomes false
if we delete the assumption that $F_\V(x_1,\dots,x_m)$ generates $\Pv.$
In the absence of that assumption, it is not clear what forms
the relation of mutual embeddability can assume.

\begin{question}\label{Q.eqrel}
For $\V$ a variety, let us call natural numbers $m$ and $n$
$\!\V\!$-equivalent \textup{(}with respect to embeddability
of free algebras\textup{)} if $F_\V(x_1,\dots,x_m)$ and
$F_\V(x_1,\dots,x_n)$ each contain an isomorphic copy of the other.
Clearly, the $\!\V\!$-equivalence classes are blocks of
consecutive integers.
Which decompositions of the natural numbers into blocks
can be realized in this way?

More generally, given a prevariety $\Pv$ and algebras
$A_1,\dots,A_r$ in $\Pv,$ let us define a preorder $\preceq$
on \mbox{$\!r\!$-tuples} of natural numbers by writing
$(m_1,\dots,m_r)\preceq (n_1,\dots,n_r)$ if the coproduct in $\Pv$
of $m_1$ copies of $A_1,$ $m_2$ copies of $A_2,$ etc., through
$m_r$ copies of $A_r,$ is embeddable in the coproduct
of $n_1$ copies of $A_1$ etc., through $n_r$ copies of $A_r.$
What preorderings on $\omega^r$ can be realized in this way?
In particular, what equivalence relations on $\omega^r$
can be the equivalence relation determined by such a preorder?
Do these answers change if one requires $\Pv$ to be a variety?
\end{question}

In~\cite{MVZ}, for certain varieties $\V$ of Lie algebras over
a field of characteristic $0,$ bounds are
obtained on $n/m$ for any $m,\,n$ equivalent under
the relation of the first paragraph of Question~\ref{Q.eqrel} above.
The idea is to note that if $F_\V(x_1,\dots,x_m)$ and
$F_\V(x_1,\dots,x_n)$ are mutually embeddable, they must
have the same Gelfand-Kirillov dimension (a measure of growth rate).
Upper and lower bounds are obtained for the Gelfand-Kirillov
dimension of $F_\V(x_1,\dots,x_n),$ leading to
the asserted conclusions.
However, it seems most likely that for such $\V,$ the
Gelfand-Kirillov dimension of $F_\V(x_1,\dots,x_n)$ will grow
with $n$ in a ``smooth'' fashion; if so, one should
in fact be able to prove that no free algebras of
distinct ranks in $\V$ are mutually embeddable, in which case such
varieties will not give interesting examples relevant
to Question~\ref{Q.eqrel}.

For results on {\em isomorphisms} and {\em surjections} among
free algebras, rather than embeddings, see~\cite{SS} and~\cite{IBN}.
The latter shows that all consistent cases are
realized by module-varieties $\mathbf{Mod}_R$ for rings $R.$

In generalizing Proposition~\ref{P.2->3} from free algebras to general
coproducts, we found that the context that made the
argument work was that of coproducts in a {\em prevariety}.
Theorem~\ref{T.P-indp} does not give us the corresponding statement
for general $A_0,\dots,A_n,\,B_1,\dots,B_n$
with the prevariety $\Pv$ replaced by a variety $\V,$ and
the condition that $\Pv$ be generated by $A_n$ as a prevariety replaced
by the condition that $\V$ be generated by $A_n$ as a variety,
since for a variety $\V$ and an algebra $A\in\V,$
the condition that $\V$ be generated
by $A$ as a variety is weaker than the condition that
it be generated by $A$ as a prevariety.
(E.g., the variety of abelian groups is generated by the
infinite cyclic as a variety, but not as a prevariety, since all
groups in the prevariety it generates must be torsion-free.)
But I don't have a counterexample to the modified statement.

\begin{question}\label{Q.variety}
\textup{(i)} \ Does Theorem~\ref{T.P-indp} remain true if
``prevariety'' is everywhere replaced by ``variety''?
\vspace{.2em}

If not, or if the question proves difficult,
one might examine some special cases, such as\\[.2em]
\textup{(ii)} \ If in Corollary~\ref{C.A_i=A}
we replace ``prevariety'' by ``variety'',
and add the assumption that $A$ is free of rank~$1$ in that variety
\textup{(}but not that the $B_i$ are free\textup{),}
does the statement still hold?\\[.2em]
\textup{(iii)} \ If $\V$ is a variety, and $A$ an algebra that
generates~$\V$ as a variety, and that contains as a
subalgebra a coproduct of two copies of itself in~$\V,$ must
it contain a coproduct of three copies of itself in~$\V$?
\end{question}

More likely to have positive answers, since quasivarieties
are more like prevarieties than varieties are, is

\begin{question}\label{Q.qv}
Same questions \textup{(i), (ii), (iii)} as above, but with
``variety'' everywhere replaced by ``quasivariety''
\textup{(}necessarily, of finitary algebras\textup{)}.
\end{question}

Looking back further, to \S\ref{S.free},
the mixture of positive and negative results there suggests

\begin{question}\label{Q.algorithm}
Is there a nice criterion for whether
a $\!3\!$-tuple $(a,\<b,\<c)$ of monoid words in two letters
$p,\ q$ has the property proved in Proposition~\ref{P.2->3}
to hold for the $\!3\!$-tuple $(\<p,\,pq,\,pqq),$
and in Proposition~\ref{P.not_px,pqx,qqx} {\em not} to hold
for the $\!3\!$-tuple $(\<p,\,pq,\,qq),$ namely, of witnessing
the existence of subalgebras free on $3$ generators
in all relatively free $\!1\!$-generator algebras
$\langle x\rangle$ that contain $\!2\!$-generator subalgebras
$\langle\<p\<x,q\<x\rangle$ in the varieties they generate?

More generally, given $n>1$ and $N>1,$ one may ask
for a criterion for an $\!N\!$-tuple $(a_1,\dots,a_N)$
of words in $n$ letters $p_1,\dots,p_n$ to witness the existence of a
free subalgebra on $N$ generators in any relatively free algebra on one
generator that contains a free subalgebra
$\langle\<p_1x,\dots,p_nx\rangle$ on $n$ generators in the
variety it generates.

\textup{(}Still more generally, for $n>m>0$ and $N>1,$ one may ask
how to decide whether a given $\!N\!$-tuple of terms in $m$
variables and $n$ operation symbols each of arity $m$ witnesses
the result of Corollary~\ref{C.m<n<N}.
But since terms in operation symbols of arity $>1$ are more complicated
than words in unary operation symbols, there seems to be less
likelihood of a simple answer.\textup{)}
\end{question}

The next four sections are related to, but
do not depend on, the material above, except for the definitions.
Section~\ref{S.amalg} recalls what it means for a category
of algebras to have the {\em amalgamation property}, obtains some
equivalent statements, and then shows that for prevarieties
with that property, one has stronger results on when a family
of subalgebras of an algebra generates a subalgebra isomorphic to
their coproduct than those that we have seen above to hold in general.
In a different direction, motivated by the fact
that the prevarieties considered in \S\ref{S.P-indp}
were by hypothesis each generated by a single algebra,
\S\S\ref{S.compat1}-\ref{S.comfort} show that the
number of algebras needed to generate a prevariety has important
consequences for the behavior of coproducts therein.
The brief section \S\ref{S.Sym}, which is included in this note
only for convenience, answers a different question about coproducts,
also raised in~\cite{embed}, concerning subgroups of
the full symmetric group on an infinite set.

\section{The amalgamation property, and its consequences for $\!\Pv\!$-independence.}\label{S.amalg}

In any class of algebras that admits coproducts with amalgamation
(pushouts), it is well known and easy to verify
that the {\em amalgamation property}
(definition recalled in~(\ref{d.amalg}) below)
is equivalent to the condition that for all pairs of one-to-one
maps with common domain, $A\to B$ and $A\to C,$ the coprojection
maps of $B$ and $C$ into the coproduct with amalgamation
$B\,\cP_A\,C$ are also one-to-one.
The next lemma gives some further consequences of that
property, in the same vein.
We formulate it in a context more general than that of
categories of algebras, though less sophisticated than
that \mbox{of~\cite[\S6]{KMPT}}.

In that lemma, the functor $U:\mathbf{C}\to\mathbf{Set}$
plays the role of the underlying set functor of a category of
algebras, but we shall not need to assume it faithful,
as one does when defining the concept of a concrete category.

One other notational remark: So far, I have generally written
$\coprod_\Pv^{}$ for ``coproduct in the category $\Pv$''; but when
discussing coproducts with amalgamation of an object, we
will use the subscript position for that object, leaving the
category to be understood from the context.
I will follow this mixed practice for the rest of the paper,
showing the category when no amalgamation is involved.
(The superscript position, which might otherwise be assigned to one
of these, is used here for index sets over which coproducts are taken.
If there were danger of ambiguity, we could write
$\coprod_{\Pv,\<A}^{}$ rather than $\coprod_A^{},$
or regard coproducts with amalgamation as coproducts
in a comma category $(A\downarrow\Pv)$ and so write
$\coprod_{(A\,\downarrow\,\Pv)}^{}.)$

\begin{lemma}\label{L.amalg}
Let $\mathbf{C}$ be a category and $U:\mathbf{C}\to\mathbf{Set}$
a functor, and let us call a morphism $f$ in $\mathbf{C}$
``one-to-one'' if $U(f)$ is a one-to-one set map, and emphasize
this by indicating such morphisms using tailed arrows:~$\tto.$

Assume that $\mathbf{C}$ admits pushouts of pairs of one-to-one
morphisms; i.e., that if $f:S\tto A$ and $g:S\tto B$ are one-to-one,
then the coproduct with amalgamation $A\,\cP_S\,B$ exists.
\textup{(}But we do not assume at this point that the maps of $S,$
$A$ and $B$ to that coproduct are one-to-one.\textup{)}

Then the following three conditions are equivalent:
\begin{equation}\begin{minipage}[c]{35pc}\label{d.amalg}
$\mathbf{C}$ has the amalgamation property
\textup{(\cite[p.82]{KMPT})}.
That is, given objects $A,\ B,\ C$ of $\mathbf{C},$ and
one-to-one morphisms $f:A\tto B,\ g:A\tto C,$ there exists
an object $D,$ and one-to-one morphisms $f':B\tto D,\ g':C\tto D,$
such that $f'f=g'\,g.$
\end{minipage}\end{equation}
\begin{equation}\begin{minipage}[c]{35pc}\label{d.SAB,ST}
For all objects $S,\,T,\,A,\,B$ and one-to-one morphisms
$S\tto T,$ $S\tto A,$ and $f:A\tto B$ in $\mathbf{C},$ the
induced morphism $f\cP_S T:\,A\cP_S\,T\to B\cP_S\,T$ is one-to-one.
\end{minipage}\end{equation}
\begin{equation}\begin{minipage}[c]{35pc}\label{d.nSAiBi}
For all objects $S,$ positive integers $n,$ and finite families of
objects and one-to-one morphisms $S\tto A_i$ and $f_i:A_i\tto B_i$
in $\mathbf{C}$ $(i\<{=}\<1,\dots,n),$ the induced morphism
$\coprod_S^{i=1,\dots,n}f_i:\,\coprod_S^{i=1,\dots,n}A_i\,\to\,
\coprod_S^{i=1,\dots,n}B_i$ is one-to-one.
\end{minipage}\end{equation}

Moreover, if $\mathbf{C}$ also admits direct limits \textup{(}colimits
over directed partially ordered sets\textup{)}, and if $U$ respects
these \textup{(}e.g., if $\mathbf{C}$ is a
quasivariety of finitary algebras, and $U$ its
underlying set functor\textup{)}, then $\mathbf{C}$ has coproducts
with amalgamation of possibly infinite families
of one-to-one maps $S\tto A_i$ \textup{(}for fixed $S,$ and $i$ ranging
over a possibly infinite set $I);$ and~\textup{~(\ref{d.nSAiBi})}
goes over to such coproducts.
That is,~\textup{(\ref{d.amalg})-(\ref{d.nSAiBi})}
are also equivalent to:
\begin{equation}\begin{minipage}[c]{35pc}\label{d.ISAiBi}
For all objects $S,$ nonempty sets $I,$ and families of
objects and one-to-one morphisms $S\tto A_i$ and $f_i:A_i\tto B_i$ in
$\mathbf{C}$ $(i\in I),$ the induced morphism $\coprod_S^{i\in I}f_i:\,
\coprod_S^{i\in I}A_i\to\coprod_S^{i\in I}B_i$ is one-to-one.
\end{minipage}\end{equation}
\end{lemma}

\begin{proof}
(\ref{d.amalg})$\Rightarrow$(\ref{d.SAB,ST}):
Given objects and maps as in~(\ref{d.SAB,ST}), the
amalgamation property implies, as mentioned above, that the
coprojection $A\to A\cP_S\,T$ is one-to-one.
From this and the assumed one-to-one-ness of the map $A\to B$
we similarly get one-to-one-ness of the coprojection
$A\cP_S T\to B\cP_A(A\cP_S T)=B\cP_S T,$ as desired.

(\ref{d.SAB,ST})$\Rightarrow$(\ref{d.amalg}): Given objects and
maps as in~(\ref{d.amalg}), apply~(\ref{d.SAB,ST})
with $A$ and its identity map in the role of $S$ and its map to $A,$
and with $C$ in the role of $T,$ noting that the domain of the resulting
map, $A\cP_A C,$ can be identified with $C.$
This gives one-one-ness of the coprojection $C\to B\cP_A\,C.$
By symmetry one also has one-one-ness of the coprojection
$B\to B\cP_A\,C.$
Taking $D=B\cP_A C,$ we get~(\ref{d.amalg}).

(\ref{d.SAB,ST})$\Rightarrow$(\ref{d.nSAiBi}):
The case $n=1$ of~(\ref{d.nSAiBi}) is trivial.
To get the case $n=2$ we make a double application of~(\ref{d.SAB,ST}),
first getting one-one-ness
for $f_1\cP_S\,A_2: A_1\cP_S\,A_2\tto B_1\cP_S\,A_2$ and then
for $B_1\cP_S\,f_2: B_1\cP_S\,A_2\tto B_1\cP_S\,B_2.$
Composing, we get one-one-ness of the desired map.

This shows that two-fold coproducts over $S$ respect one-to-one-ness
of maps among objects having one-to-one maps of $S$ into them.
Induction now gives
the corresponding result for $\!n\!$-fold coproducts.

(\ref{d.nSAiBi})$\Rightarrow$(\ref{d.SAB,ST}):
Given objects and maps as in~(\ref{d.SAB,ST}), apply the $n=2$ case
of~(\ref{d.nSAiBi}) with $A\tto B$ in the role of $f_1:A_1\tto B_1$
and the identity map of $C$ in the role of $f_2:A_2\tto B_2.$

Under the additional assumptions about direct limits, one notes that
for infinite $I,$ one can obtain $\coprod_S^{i\in I}A_i$ as the direct
limit, over the directed system of all finite subsets $I_0\subseteq I,$
of the objects $\coprod_S^{i\in I_0}A_i.$
Since by~(\ref{d.nSAiBi}), the indicated maps among these
finite coproducts are one-to-one, and by assumption direct limits
respect $U$ (and hence one-one-ness), the corresponding maps among
the coproducts over $I$ are also one-to-one.
The converse is immediate:~(\ref{d.ISAiBi}) includes~(\ref{d.nSAiBi}).
\end{proof}

Let us now note how the amalgamation property implies conditions on
independent subalgebras stronger than those of~\S\ref{S.P-indp}.
In considering categories of algebras, we shall understand the
functor $U$ of Lemma~\ref{L.amalg} to be the underlying set functor.
Thus, ``one-to-one'', in our formulation~(\ref{d.amalg}) of the
amalgamation property and our statements of conditions equivalent
thereto, has its usual meaning for algebras.

\begin{corollary}\label{C.cPincP}
Suppose that $\Pv$ is a prevariety having the amalgamation property,
that $A$ is a $\!\Pv\!$-algebra, that $(B_i)_{i\in I}$ is a finite
$\!\Pv\!$-independent family of subalgebras of $A,$ and that for
each $i\in I,$ $(C_{ij})_{j\in J_i}$ is a finite
$\!\Pv\!$-independent family of subalgebras of $B_i.$
Then $(C_{ij})_{i\in I,\,j\in J_i}$ is a
$\!\Pv\!$-independent family of subalgebras of~$A.$
\textup{(}In particular, in such a prevariety, examples like
those of Propositions~\ref{P.2_not_3}
and~\ref{P.not_px,pqx,qqx} cannot occur.\textup{)}

If $\Pv$ is in fact a {\em quasivariety} having the amalgamation
property, then the above result
holds without the finiteness restrictions on $I$ and the~$J_i.$
\end{corollary}

\begin{proof}
Since all the algebras named are subalgebras of $A,$
the unique homomorphic image of the initial algebra of $\Pv$
in all of them is the same; let us call this $S.$
Because $S$ is a homomorphic image of the initial algebra of our
category, the operator $\coprod_S$ on nonempty families
of algebras containing $S$ is just $\coprod_\Pv.$

We now apply the implication
(\ref{d.amalg})$\Rightarrow$(\ref{d.nSAiBi}) of Lemma~\ref{L.amalg}
(or if $\Pv$ is a quasivariety, the stronger implication
(\ref{d.amalg})$\Rightarrow$(\ref{d.ISAiBi})), taking
for the $S$ of~(\ref{d.nSAiBi}) and~(\ref{d.ISAiBi})
the $S$ of the preceding paragraph, and for the maps $A_i\tto B_i$
the inclusions $\coprod_\Pv^{j\in J_i}C_{ij}\subseteq B_i.$
We conclude that the natural map
from $\coprod_\Pv^{i\in I}(\coprod_\Pv^{j\in J_i}C_{ij})=
\coprod_\Pv^{i\in I,\,j\in J_i}C_{ij}$ to $\coprod_\Pv^I B_i$
is one-to-one.
Identifying the latter algebra with its embedded
image in $A,$ we get the desired conclusion.

To get the parenthetical remark about examples like those
of Propositions~\ref{P.2_not_3} and~\ref{P.not_px,pqx,qqx},
we take $I=\{0,1\},\ J_0=\{0\}, J_1=\{0,1\},$ and let
$A,$ the $B_i$ and the $C_{ij}$ all be free of rank~$1.$
\end{proof}

Remarks: Lemma~\ref{L.amalg} was a compromise between
proving the minimum we needed to get the above corollary --
that~(\ref{d.amalg}) implies the special case of~(\ref{d.nSAiBi})
where $S$ is the image of the initial object of $\mathbf{C}$ in $A,$
and so can be ignored in forming coproducts
(and if $\mathbf{C}$ has, and $U$ respects,
direct limits, the corresponding case of~(\ref{d.ISAiBi})) --
and digressing to state and prove a more complete statement.
That statement would involve the versions of conditions~(\ref{d.amalg})
and~(\ref{d.nSAiBi}) for $\!\kappa\!$-fold families for any
cardinal $\kappa,$ would establish the equivalence between those
two conditions for each such
$\kappa,$ would note that the statements for larger
$\kappa$ imply those for smaller $\kappa,$ and would verify that
the statements for {\em finite} $\kappa\geq 2$ are all equivalent, and
also equivalent to~(\ref{d.SAB,ST}).
The reader should not find it hard to work out the details.

The reason we brought $S$ into~(\ref{d.nSAiBi}), though the only case
of~(\ref{d.nSAiBi}) that our application needed was where $S$ was
a homomorphic image of
the initial object and so had no effect, was so as to get an
if-and-only-if relation between~(\ref{d.nSAiBi}) and~(\ref{d.amalg}),
the amalgamation property.
(The latter is a well-known property, satisfied
by the categories of groups, semilattices, lattices, and
commutative integral domains, and many others.
See the first column of the table in
\cite[pp.98-107]{KMPT} for more results, positive and negative.)
That equivalence fails if $S$ in~(\ref{d.nSAiBi}) is
restricted to homomorphic images
of the initial object.
For instance, the normal form for coproducts of monoids
shows that the variety $\mathbf{Monoid}$ satisfies
the cases of~(\ref{d.SAB,ST})-(\ref{d.ISAiBi})
where $S$ is the initial (trivial) monoid.
However $\mathbf{Monoid}$ does not satisfy the amalgamation
property~(\ref{d.amalg});
e.g., letting $A=\langle x\rangle,$ the free monoid on one generator,
and letting $B$ and $C$ be the overmonoids of $A$ gotten by
adjoining a left inverse $y,$ respectively a right inverse $z,$ to $x,$
one finds that in $B\<\coprod_A C,$ the elements $xy$ of $B$ and
$zx$ of $C$ fall together with $1;$ so the maps from $B$ and
$C$ to this algebra are not one-to-one.
On the other hand, because the special case
of~(\ref{d.SAB,ST})-(\ref{d.ISAiBi}) which we have seen
suffices for Corollary~\ref{C.cPincP} holds,
$\mathbf{Monoid}$ does satisfy the conclusion of that corollary.

Here is another result (alluded to in the discussion
following Lemma~\ref{L.dirsys}) of a sort similar to the above,
which for simplicity of wording we will again state in terms
of the amalgamation property, though again, only the cases
of~(\ref{d.SAB,ST})-(\ref{d.ISAiBi}) where $S$ is a homomorphic
image of the initial object of $\Pv$ are needed.

\begin{corollary}\label{C.indincind}
Suppose $\Pv$ is a prevariety having the amalgamation property,
and $A$ a $\!\Pv\!$-algebra.
Then every nonempty subfamily of a $\!\Pv\!$-independent
family of subalgebras of $A$ is $\!\Pv\!$-independent.
\end{corollary}

\begin{proof}
Given a $\!\Pv\!$-independent family of subalgebras
$B_i$ $(i\in I)$ of $A,$ their $\!\Pv\!$-independence says
that the subalgebra of $A$ that they generate is
isomorphic to their coproduct, which we see coincides with
their coproduct over the common image $S$ in all these algebras of the
initial algebra of $\Pv.$
For any nonempty subset $J\subseteq I,$ the coproduct
of the $B_i$ for $i\in J$ likewise coincides with
their coproduct over $S.$
We now apply~(\ref{d.SAB,ST})
with this algebra $S$ for both the $S$ and $A$ of that condition,
with $\coprod_S^J B_i$ for $T,$ and with $\coprod_S^{I-J} B_i$ for $B,$
and then bring in the assumed $\!\Pv\!$-independence of the whole
family.
We thus get one-one-ness of the natural maps shown by the first arrow in
\begin{equation}\begin{minipage}[c]{35pc}\label{d.JtoI}
$\coprod_\Pv^J B_i\ \cong\ \coprod_S^J B_i\ \tto
\ \coprod_S^J B_i\cP_S\coprod_S^{I-J}\,B_i\ \cong
\ \coprod_S^I B_i\ \cong\ \coprod_\Pv^I B_i\ \tto\ A,$
\end{minipage}\end{equation}
and the above arrows and
isomorphisms compose to the map we wished to show one-to-one.
\end{proof}

\section{$\!\Pv\!$-compatible algebras.}\label{S.compat1}

The prevarieties considered in \S\ref{S.P-indp} were
each generated by a single algebra.
Although any {\em variety} of algebras can be generated
as a {\em variety} by a single algebra (namely, by a free algebra on
sufficiently many generators), prevarieties generated
as {\em prevarieties} by a single algebra are rather special.
This was shown by Mal'cev for quasivarieties,
in a result that we will generalize in the next section.
In this section we shall see that the size of the collection of algebras
needed to generate $\Pv$ as a prevariety is a nontrivial and
interesting invariant of $\Pv,$ even if $\Pv$ happens to be a variety.

\begin{definition}\label{D.compat}
Let $\Pv$ be a prevariety.
Then a set $\X$ of $\!\Pv\!$-algebras will be called
{\em $\!\Pv\!$-compatible} if for every $A_0\in\X,$
the coprojection map $A_0\to\coprod_\Pv^{A\in\X} A$ is one-to-one;
equivalently, if there exists an algebra $B$ in $\Pv$
admitting one-to-one homomorphisms $A\to B$ for all $A\in\X.$
\end{definition}

\begin{theorem}\label{T.compat1}
Suppose $\Pv$ is a prevariety that is
residually small \textup{(}i.e., that can generated as a prevariety by
a {\em set} of algebras\textup{)} and $\kappa$ is a cardinal.
Then condition\textup{~(\ref{d.*kgen})} below implies
condition\textup{~(\ref{d.*kdecomp})}; and if $\Pv$ is a quasivariety
\textup{(}in which case, we recall, our algebras are assumed
finitary\textup{)}, the two conditions are equivalent.
\begin{equation}\begin{minipage}[c]{35pc}\label{d.*kgen}
$\Pv$ can be generated, as a
prevariety, by a set of $\leq\kappa$ algebras.
\end{minipage}\end{equation}
\begin{equation}\begin{minipage}[c]{35pc}\label{d.*kdecomp}
Every set $\X$ of subdirectly irreducible algebras
in $\Pv$ can be written as the union of $\leq\kappa$ subsets
$\X_\alpha$ $(\alpha\in\kappa),$ each of which is $\!\Pv\!$-compatible.
\end{minipage}\end{equation}
\end{theorem}

\begin{proof}
(\ref{d.*kgen})$\Rightarrow$(\ref{d.*kdecomp}):
Suppose $\Pv$ is generated by a set of $\leq\kappa$ algebras,
$\mathbf{Y}=\{B_\alpha\mid\alpha\in\kappa\},$
and that as in~(\ref{d.*kdecomp}),
$\X$ is a set of subdirectly irreducible algebras in $\Pv.$
Each $A\in\X$ is embeddable in a direct product of copies of
the $B_\alpha,$ hence, being subdirectly irreducible,
in one of the $B_\alpha.$
Letting $\X_\alpha$ be the set of members of $\X$ embeddable in
$B_\alpha,$ we get the conclusion of~(\ref{d.*kdecomp}) (using the
second formulation in the definition of $\!\Pv\!$-compatibility).

To prove that when $\Pv$ is a quasivariety,
(\ref{d.*kdecomp})$\Rightarrow$(\ref{d.*kgen}), note that by our
residual smallness hypothesis,
there is a set $\X$ of subdirectly irreducible
algebras in $\Pv$ which contains, up to isomorphism, all such algebras.
By~(\ref{d.*kdecomp}) we may write $\X=
\bigcup_{\alpha\in\kappa}\X_\alpha$
where each $\X_\alpha$ is $\!\Pv\!$-compatible.
Hence for each $\alpha,$
we can choose an algebra $A_\alpha$ in $\Pv$ in which all members of
$\X_\alpha$ can be embedded.
Since a quasivariety is generated as a prevariety by its subdirectly
irreducible algebras \cite[Theorem~3.1.1]{VAG}, the prevariety
generated by $\{A_\alpha\mid\alpha\in\kappa\}$ is all of~$\Pv.$
\end{proof}

To get easy examples showing that the least $\kappa$
for which~(\ref{d.*kdecomp}) holds
can be, inter alia, any natural number,
consider algebras with a single unary operation
$a,$ and for each positive integer $d,$ let $C_d$ be the algebra
of this type consisting of $d$ elements,
$x,\,ax,\,\dots,\,a^{d-1}x,$ cyclically permuted by $a.$

Now let $n$ be any natural number,
and $d_1,,\dots,d_n$ be positive integers none of which
is the least common multiple of any subset of the others.
(In particular, none of them is $1,$ since $1$ is the
least common multiple of the empty set.)
Let $\Pv$ be the prevariety generated by the $n$
algebras $C_{d_1},\dots,C_{d_n}.$
Since this is generated by finitely many finite finitary algebras,
it is a quasivariety.
From the description $\Pv=\Su\Pr\{C_{d_1},\dots,C_{d_n}\}$
we see that all algebras in $\Pv$ satisfy
\begin{equation}\begin{minipage}[c]{35pc}\label{d.a^d}
$(\forall\,x)\ \ a^{\mathrm{lcm}(d_1,\dots,d_n)}x\ =\ x,$
\end{minipage}\vspace{-0.4em}\end{equation}
\begin{equation}\begin{minipage}[c]{35pc}\label{d.ax=x}
$(\forall\,x,\,y,\,z)\ \ ax\ =\ x\implies y\ =\ z,$
\end{minipage}\vspace{-0.4em}\end{equation}
\begin{equation}\begin{minipage}[c]{35pc}\label{d.lcm}
For all $x,$ the least positive integer $d$ such that $a^dx=x$ is\\
the least common multiple of some subset of~$\{d_1,\dots,d_n\},$
\end{minipage}\vspace{-0.4em}\end{equation}
\begin{equation}\begin{minipage}[c]{35pc}\label{d.only1}
$(\forall\,x,y)\ \ a^dx=x\implies a^dy=y.$
\end{minipage}\end{equation}
From~(\ref{d.a^d})-(\ref{d.only1}) and our assumption that none of
the $d_i$ is the least common multiple of a subset of the rest, one
can verify that the subdirectly irreducible objects of $\Pv$ are
precisely the $n$ algebras $C_{d_i};$ and by~(\ref{d.only1})
these are pairwise incompatible; so for this
quasivariety, the least $\kappa$ as in Theorem~\ref{T.compat1} is $n.$

For $d_1,\dots,d_n$ as above, consider next, for contrast, the
quasivariety $\Pv$ generated by a single algebra,
the disjoint union $C_{d_1}\sqcup\,\dots\,\sqcup\,C_{d_n}.$
This will still satisfy~(\ref{d.a^d}), (\ref{d.ax=x}) and~(\ref{d.lcm}),
but not~(\ref{d.only1}).
The algebras $C_{d_1},\dots,C_{d_n}$ will still be subdirectly
irreducible in $\Pv,$ but they are no longer incompatible.
Indeed, since $\Pv$ is generated by a single algebra, the
least cardinal $\kappa$ as in Theorem~\ref{T.compat1} is now $1.$

For an intermediate case, given $d_1,\,d_2,\,d_3$ as above, let
$\Pv$ be generated by the three disjoint unions $C_{d_1}\sqcup C_{d_2},$
$C_{d_1}\sqcup C_{d_3}$ and $C_{d_2}\sqcup C_{d_3}.$
Since none of these generating algebras contains copies
of all three $C_{d_i},$ these algebras, and hence all algebras
in $\Pv,$ satisfy the implication
\begin{equation}\begin{minipage}[c]{35pc}\label{d.not3}
$(\forall\ x_1,\,x_2,\,x_3,\,y,\,z)
\ \ (a^{d_1}x_1=x_1)\wedge(a^{d_2}x_2=x_2)
\wedge(a^{d_3}x_3=x_3)\implies y=z.$
\end{minipage}\end{equation}
Hence, though any two of $C_{d_1},$ $C_{d_2}$ and $C_{d_3}$ are
$\!\Pv\!$-compatible, the set consisting of all three is not.

If we take an infinite sequence of integers $d_1,\,d_2,\,\dots\,,$ none
of which divides any of the others (for instance, the primes), and let
$\Pv$ be the prevariety generated by all {\em finite} disjoint unions
of the $C_{d_i},$ this will no longer be a quasivariety.
For it will satisfy the sentence
\begin{equation}\begin{minipage}[c]{35pc}\label{d.inf_sent}
$(\forall\,x_1,\dots,x_n,\dots;\ y;\ z)
\ \ (\,\bigwedge_{i=1}^\infty\ a^{d_i}x_i=x_i)\implies y=z,$
\end{minipage}\end{equation}
so the direct limit, as sets with one unary operation, of the above
generating family of finite unions of $C_i$ (mapped into one
another by inclusion) does not lie in $\Pv;$ their direct limit
in $\Pv$ is, by~(\ref{d.inf_sent}), the trivial (one-element) algebra.
This example also shows that in our earlier result, Lemma~\ref{L.amalg},
the added direct limit hypothesis was indeed needed
to get from (\ref{d.amalg})-(\ref{d.nSAiBi}) to~(\ref{d.ISAiBi}).
For it is easy to see that $\Pv$ satisfies~(\ref{d.amalg}), while to
see that it does not satisfy~(\ref{d.ISAiBi}), we may
take for $S$ and the $A_i$ the empty algebra,
and for the $B_i$ the above algebras $C_{d_i}.$

To see that the
implication (\ref{d.*kdecomp})$\Rightarrow$(\ref{d.*kgen}) of
Theorem~\ref{T.compat1}, which we proved for
quasivarieties, does not hold for general prevarieties, let us construct
a prevariety not having ``enough'' subdirectly irreducible algebras:
Let $p$ be a prime, let $G$ be the additive group of an
infinite-dimensional vector space over the field of $p$ elements, and
let $\Pv$ be the prevariety consisting of
all $\!G\!$-sets $A$ such that if an element
of $A$ is fixed by an element of $G,$ then all elements of $A$ are
fixed by that element, and if an element of $A$ is fixed by
infinitely many elements of $G,$ then all elements of $A$ are equal.
Then $\Pv$ is residually small: the set of $\!G\!$-sets $G/H$ for
all finite subgroups $H\subset G$ generates $\Pv.$
But any finite-dimensional subspace of $G$ is an intersection of
two properly larger finite dimensional subspaces, hence any nontrivial
algebra in $\Pv$ can be decomposed as a subdirect product of algebras
with larger pointwise stabilizers; so
$\Pv$ has no subdirectly irreducible algebras, so it
satisfies~(\ref{d.*kdecomp}) for $\kappa=0.$
On the other hand, no two nonempty algebras in $\Pv$
having different pointwise stabilizers are
compatible, so~(\ref{d.*kgen}) does not hold for any finite $\kappa.$

\section{All under one roof: prevarieties where all algebras are $\!\Pv\!$-compatible.}\label{S.compat2}

For prevarieties that can be generated by {\em one} algebra, a
stronger result can be proved than the $\kappa=1$ case
of~(\ref{d.*kdecomp}); moreover, we can weaken the
above assumption ``generated by one algebra'' to
a condition that is necessary as well as sufficient
for our strengthened conclusion.

We need the following definition.
(Recall that a preordering $\preceq$ on a set means a reflexive,
transitive, but not necessarily antisymmetric binary relation.)

\begin{definition}\label{D.absdird}
A preordered {\em class} $(\mathbf{K},\preceq)$ will be called
{\em absolutely directed} if every {\em set} of elements
of $\mathbf{K}$ is majorized by some element of $\mathbf{K}.$
\end{definition}

In particular, a preordered {\em set} is absolutely
directed if and only if it has a greatest element (an
element $\succeq$ all elements).

In the next result, condition~(\ref{d.setcompat}) can
be seen to be a strengthening of the $\kappa=1$ case
of~(\ref{d.*kdecomp}) (with ``nontrivial'' replacing
``subdirectly irreducible''), while~(\ref{d.absdirgen}) is
a weakening of the condition that
$\Pv$ be generated as a prevariety by a single algebra.
The equivalence of~(\ref{d.2compat})
and~(\ref{d.qv1gen}) for quasivarieties is due to Mal'cev.

By the {\em trivial} algebra we will
always mean the one-element algebra.
(So trivial algebras and empty algebras are never the same thing.)

\begin{theorem}[{cf.\ Mal'cev \cite{AIM},
\cite[Proposition~2.1.19]{VAG}}]\label{T.absdird}
Let $\Pv$ be a prevariety, and let $\preceq$ be the preordering
``is embeddable in'' among algebras in $\Pv.$
Then the following conditions are equivalent.
\begin{equation}\begin{minipage}[c]{35pc}\label{d.setcompat}
Every set of nontrivial algebras in $\Pv$ is $\!\Pv\!$-compatible.
\end{minipage}\end{equation}
\begin{equation}\begin{minipage}[c]{35pc}\label{d.absdirgen}
$\Pv$ is generated as a prevariety by a class of
algebras absolutely directed under ${\preceq}\,.$
\end{minipage}\end{equation}

Moreover, if $\Pv$ is a quasivariety
\textup{(}so that, again, our algebras are assumed finitary\textup{)},
then the above conditions are also equivalent to each of
\begin{equation}\begin{minipage}[c]{35pc}\label{d.2compat}
Every pair of nontrivial algebras in $\Pv$ is $\!\Pv\!$-compatible.
\end{minipage}\end{equation}
\begin{equation}\begin{minipage}[c]{35pc}\label{d.qvabsdirgen}
$\Pv$ is generated as a {\em quasivariety} by a
class of algebras absolutely directed under ${\preceq}\,.$
\end{minipage}\end{equation}
\begin{equation}\begin{minipage}[c]{35pc}\label{d.qv1gen}
$\Pv$ is generated as a quasivariety by a single algebra.
\end{minipage}\end{equation}
\end{theorem}

\begin{proof}
(\ref{d.setcompat}) says that the class of nontrivial
algebras in $\Pv$ is absolutely directed under ${\preceq}\,.$
But $\Pv$ is generated as a prevariety by its nontrivial algebras
(the trivial algebra being the direct product of the empty family
thereof), so this implies~(\ref{d.absdirgen}).

To show (\ref{d.absdirgen})$\Rightarrow$(\ref{d.setcompat}),
let $\X$ be an absolutely directed class of algebras generating
$\Pv,$ and $\mathbf{Y}$ any set of nontrivial algebras in $\Pv.$
Since $\mathbf{Y}$ is a set, we can find a set $\X_0\subseteq\X,$
homomorphisms into members of which separate points of algebras
in $\mathbf{Y},$ and by the
absolute directedness of $\X,$ some one algebra $A\in\X$
contains embedded images of all members of $\X_0;$ hence homomorphisms
into $A$ separate points of algebras in $\mathbf{Y}.$
Hence if we form a direct product $A^I$ of sufficiently many copies of
$A,$ then for each
nonempty $B\in\mathbf{Y},$ we can use maps to some coordinates
of $A^I$ to separate points of $B;$ and since $B$ is nontrivial and
nonempty, the set of maps so used will be nonempty, and we can
repeat some of them to fill in the remaining coordinates if any; thus
we can embed $B$ in $A^I.$
The same conclusion is vacuously true if $B$ is empty,
so $A^I$ has subalgebras isomorphic to all $B\in\mathbf{Y},$
proving that $\mathbf{Y}$ is $\!\Pv\!$-compatible.

Now let $\Pv$ be a quasivariety.

Clearly,~(\ref{d.setcompat})$\Rightarrow$(\ref{d.2compat}).
The converse holds because we can go from pairwise coproducts
to finite coproducts by induction, while coproducts of infinite families
are direct limits of coproducts of finite families, and in
a quasivariety, direct limits respect the underlying set functor.
Thus,~(\ref{d.setcompat})-(\ref{d.2compat}) are equivalent.

(\ref{d.absdirgen})$\Rightarrow$(\ref{d.qvabsdirgen}) is
trivial, since the quasivariety generated by a class of algebras
contains the prevariety generated by the same class.
We shall now show (\ref{d.qvabsdirgen})$\Rightarrow$(\ref{d.qv1gen}),
then note two alternative ways of getting back from~(\ref{d.qv1gen}):
to~(\ref{d.absdirgen}), or to~(\ref{d.2compat}).

Given~(\ref{d.qvabsdirgen}), let $\X$ be an absolutely directed
class of algebras generating~$\Pv$ as a quasivariety.
Since the {\em finite} sentences~(\ref{d.Horn}) form (modulo
notation) a {\em set}, if we choose for each such sentence not
satisfied by~$\Pv$ a member of $\X$ for which it fails,
we get a set of algebras $\X_0\subseteq\X$ which again generates $\Pv.$
By assumption, $\X_0$ is majorized by an algebra $A\in\X,$
and this will likewise generate $\Pv,$ proving~(\ref{d.qv1gen}).

Now assume~(\ref{d.qv1gen}), and let $A$
be an algebra that generates $\Pv$ as a quasivariety.

On the one hand, one can deduce~(\ref{d.absdirgen}) from the
fact that $\Pv$ is generated as a prevariety by
the class of ultrapowers of $A$ (\cite[Corollary~2.3.4(i)]{VAG};
cf.\ last paragraph of Definition~\ref{D.qv} above) by
verifying that that class is absolutely directed under ${\preceq}\,.$
The idea is that given a set of ultrafilters $\mathcal{U}_j$ $(j\in J),$
each on a set $I_j,$ these yield a ``product'' filter on
$\prod^J I_j,$ and any ultrafilter $\mathcal{U}$ containing this
will have the property that all the ultrapowers $A^{\mathcal{U}_j}$
embed in the ultrapower $A^\mathcal{U}.$

To get~(\ref{d.2compat}), on the other hand, suppose by way of
contradiction that $B_0$ and $B_1$ were a non-$\!\Pv\!$-compatible
pair of nontrivial algebras in $\Pv.$
Without loss of generality, suppose $B_0,$
has non-one-to-one coprojection into $B_0\cP_{\Pv}\,B_1;$
let elements $x\neq y$ of $B_0$ fall together there.
Since $\Pv$ is determined by finite sentences~(\ref{d.Horn}),
the conjunction of finitely many of these universal sentences
with finitely many equations holding among finitely many elements
of $B_0$ and $B_1$ must imply $x=y.$
But every finite system of relations among elements of each of the
$B_i$ is realizable by relations among some family of elements of
$A$ (otherwise $A,$ and hence $\Pv,$ would
satisfy an implication~(\ref{d.Horn}) saying that the conjunction
of such a system of relations implies that all elements are equal,
contradicting our assumption that the $B_i$ are nontrivial).
On the other hand, since $\Pv$ does not satisfy an implication
forcing $x=y$ to hold in $B_0,$ the above equations involving
elements of $B_0$ must be satisfiable by a family of elements
of $A$ with {\em distinct} elements representing $x$ and $y.$
But combining this family with the family of elements of $A$
chosen above to satisfy our finitely many relations holding in
$B_1,$ we see that the sentences~(\ref{d.Horn})
defining $\Pv$ imply that those two
elements {\em are} equal, giving the required contradiction.
\end{proof}

In the above theorem we had to exclude the trivial
algebra from certain statements.
The following addendum to that theorem shows that in many
prevarieties, not only is that restriction unnecessary,
but trivial algebras can be used in formulating a very simple
criterion,~(\ref{d.allcompatwtriv}), for the equivalent
conditions of the theorem to hold.

\begin{corollary}\label{C.dirdwtriv}
In the context of Theorem~\ref{T.absdird}, suppose that
\begin{equation}\begin{minipage}[c]{35pc}\label{d.trivinnontriv}
$\Pv$ has at least
one {\em nontrivial} algebra with a {\em trivial} subalgebra
\end{minipage}\end{equation}
\textup{(}i.e., a nontrivial algebra with an element idempotent under
all the algebra operations\textup{)}.

Then\textup{~(\ref{d.setcompat})}, and, if $\Pv$ is a
quasivariety,\textup{~(\ref{d.2compat})}
are equivalent to the condition obtained by deleting the word
``nontrivial'' from~\textup{(\ref{d.setcompat})}; and also to
\begin{equation}\begin{minipage}[c]{35pc}\label{d.allcompatwtriv}
Every algebra in $\Pv$ is $\!\Pv\!$-compatible with the trivial algebra.
\end{minipage}\end{equation}
\end{corollary}

\begin{proof}
It is easy to deduce from~(\ref{d.trivinnontriv}) that each
of~(\ref{d.setcompat}) and~(\ref{d.2compat}) is equivalent to the
formally strengthened version of itself gotten by
deleting the restriction ``nontrivial'':
given a set $\X$ of algebras (respectively, a pair of
algebras) of $\Pv$ including the trivial
algebra, which we want to embed simultaneously in some algebra, we
``sneak the trivial member of our set in'' by hiding it in a
nontrivial algebra as in~(\ref{d.trivinnontriv}),
then apply~(\ref{d.setcompat}) (respectively,~(\ref{d.2compat})) to
the resulting family of nontrivial algebras.
As a special case of this version
of either condition, we have~(\ref{d.allcompatwtriv}).

On the other hand, given~(\ref{d.allcompatwtriv}), we can
get the strengthened form
of~(\ref{d.setcompat}) by a version of the construction by which one
embeds a family of groups in their direct product group.
Let $\{B_i\mid i\in I\}$ be any set of algebras in $\Pv.$
By~(\ref{d.allcompatwtriv}), embed each $B_i$ in an algebra $A_i$
containing an idempotent element $e_i.$
Taking $A=\prod^I A_i,$ we can embed each $B_j$ in $A$ by using
the inclusion map at the $\!j\!$-th component, and
mapping to every other component by collapsing everything
to a trivial subalgebra.
\end{proof}

Clearly {\em every} prevariety of groups, monoids, or lattices
satisfies~(\ref{d.allcompatwtriv}), hence satisfies~(\ref{d.setcompat})
with the nontriviality condition deleted,~(\ref{d.absdirgen}),
and, if it is a quasivariety,~(\ref{d.qv1gen}).

On the other hand, the variety $\V$ of unital associative
(or unital associative commutative) algebras over any field
satisfies~(\ref{d.setcompat})
(and hence~(\ref{d.absdirgen})-(\ref{d.qv1gen})), by the standard
description of coproducts of such algebras, but
not~(\ref{d.trivinnontriv}) or the version of~(\ref{d.setcompat})
with ``nontrivial'' deleted; rather, the trivial algebra in $\V$ is not
$\!\V\!$-compatible with any other algebra.
Hence in the absence of~(\ref{d.trivinnontriv}),
the exclusion of the trivial algebra in~(\ref{d.setcompat})
and~(\ref{d.2compat}) is indeed needed to make
Theorem~\ref{T.absdird} hold.
Our constructions in the preceding section with unary algebras
also illustrate this: in the prevariety generated by a single
algebra $C_d$ $(d>1),$ the conditions of Theorem~\ref{T.absdird}
hold, but $C_d$ satisfies $(\forall\,x,\,y,\,z)
\linebreak[0]\ ax\,{=}\,x{\implies}\linebreak[0]y\,{=}\,z,$ so the
trivial algebra
is not $\!\Pv\!$-compatible with any nontrivial algebra.

There are also examples where~(\ref{d.trivinnontriv})
holds, but where the equivalent conditions of
Theorem~\ref{T.absdird} and Corollary~\ref{C.dirdwtriv}
do not; for example, the variety of groups
or monoids with one distinguished element, or of lattices with two
distinguished elements: a counterexample to~(\ref{d.allcompatwtriv})
is given by any group or monoid with distinguished element that is
not the identity, or any lattice with a pair of distinguished
elements that are not equal.

These same examples also show that the analog of the
implication (\ref{d.qv1gen})$\!\implies\!$(\ref{d.setcompat}) does not
hold for varieties, with ``generated as a variety'' in place of
``generated as a quasivariety'', since every variety
satisfies the analog of~(\ref{d.qv1gen}).\vspace{.2em}

The next corollary is a result promised in
the comment following Lemma~\ref{L.dirsys}.

\begin{corollary}\label{C.subindep}
Suppose $\Pv$ is a prevariety generated by a single algebra,
or, more generally, satisfying\textup{~(\ref{d.absdirgen})},
and let $(B_i)_{i\in I}$ be a family of nontrivial algebras in $\Pv.$
\textup{(}Again, if $\Pv$
satisfies\textup{~(\ref{d.trivinnontriv}),}
the restriction ``nontrivial'' can be dropped.\textup{)}
Then
\begin{equation}\begin{minipage}[c]{35pc}\label{d.JtoI11}
For every $J\subseteq I,$ the natural map
$\coprod_\Pv^{i\in J}B_i\to\coprod_\Pv^{i\in I}B_i$ is one-to-one.
\end{minipage}\end{equation}
Hence
\begin{equation}\begin{minipage}[c]{35pc}\label{d.subfamilyindp}
Any subfamily of a $\!\Pv\!$-independent family of
subalgebras of a nontrivial algebra $A$ is\\
$\!\Pv\!$-independent.
\end{minipage}\end{equation}
\end{corollary}

\begin{proof}
To see~(\ref{d.JtoI11}), note that $\coprod_\Pv^{i\in I}B_i\cong
(\coprod_\Pv^{i\in J}B_i)\,\cP_\Pv\,(\coprod_\Pv^{i\in I-J}B_i),$ with
the natural map $\coprod_\Pv^{i\in J}B_i\to\coprod_\Pv^{i\in I}B_i$
corresponding to the first coprojection under this decomposition.
By the implication (\ref{d.absdirgen})$\Rightarrow$(\ref{d.setcompat})
(or its modified version
if $\Pv$ satisfies~(\ref{d.trivinnontriv})),
the indicated coproducts over $J$ and $I-J$ are $\!\Pv\!$-compatible,
hence that coprojection map is one-to-one, as claimed.
(A slight hiccup in this argument: If $J$ or $I-J$ is empty,
can we be sure the coproduct over that subset, namely the
initial algebra, is nontrivial?
No, but if it is trivial, and if $\Pv$ is not the trivial prevariety,
then since the initial algebra can be mapped into every
algebra,~(\ref{d.trivinnontriv}) holds, and so we
are in the case where we don't need nontriviality.)

To get~(\ref{d.subfamilyindp}), recall that the statement that
$(B_i)_{i\in I}$ is an independent family of subalgebras
of $A$ means that the subalgebra of $A$
generated by these subalgebras can be identified with their coproduct.
If none of the $B_i$ is trivial, then this observation together
with~(\ref{d.JtoI11}) immediately gives the desired conclusion.
If at least one of the $B_i$ is trivial, then since by assumption
$A$ is not,~(\ref{d.trivinnontriv}) holds, and by the parenthetical
addendum to the first part of this corollary,
we again have~(\ref{d.JtoI11}) and can proceed as before.
\end{proof}

On a different topic, let us note
the extent to which Theorem~\ref{T.compat1} does and does
not go over from prevarieties to quasivarieties.

\begin{corollary}\label{C.qv*k}
If $\Pv$ is a {\em quasivariety}, then \textup{(}even without the
residual smallness assumption of Theorem~\ref{T.compat1}\textup{)},
condition\textup{~(\ref{d.*kdecomp})} implies
\begin{equation}\begin{minipage}[c]{35pc}\label{d.qv*k}
$\Pv$ can be generated as a quasivariety by $\leq\kappa$ algebras.
\end{minipage}\end{equation}
The reverse implication holds
if $\kappa$ is finite, but not for any infinite $\kappa.$
\end{corollary}

\begin{proof}
(\ref{d.*kdecomp})$\Rightarrow$(\ref{d.qv*k}):  A quasivariety
$\Pv$ is generated as a prevariety, and hence as a quasivariety,
by its subdirectly irreducible algebras
\cite[Theorem~3.1.1]{VAG}, hence, as in
the proof of Theorem~\ref{T.absdird}, we can find a {\em set}
$\X$ of these that generates it as a quasivariety.
By~(\ref{d.*kdecomp}) we can write $\X$ as
$\bigcup_{\alpha\in\kappa}\X_\alpha$ where each $\X_\alpha$
is $\!\Pv\!$-compatible.
If for each $\alpha\in\kappa$ we let $A_\alpha$ be an algebra in $\Pv$
containing embedded images of all members of $\X_\alpha,$ then $\Pv$
is generated as a quasivariety by this set of $\kappa$ algebras.

For the converse assertion when $\kappa$ is a natural number $n,$
let $\Pv$ be generated as a quasivariety by $A_1,\dots,A_n.$
Then $\Pv=\Su\Pr\Pr_\mathrm{\!ult}\{A_1,\dots,A_n\},$ where
$\Pr_\mathrm{\!ult}$
denotes closure under ultraproducts; thus, each $\!\Pv\!$-subdirectly
irreducible object of $\Pv$ is embeddable in a member
of $\Pr_\mathrm{\!ult}\{A_1,\dots,A_n\}.$
Moreover, the operator $\Pr_\mathrm{\!ult}$ respects finite
decompositions; that is, any ultraproduct of a family of structures
indexed by a finite union of sets, $(A_i)_{i\in I_1\cup\dots\cup I_n},$
can be written as an ultraproduct of one of the subfamilies
$(A_i)_{i\in I_m}.$
Hence $\Pr_\mathrm{\!ult}\{A_1,\dots,A_n\}=
\Pr_\mathrm{\!ult}\{A_1\}\cup\dots\cup\Pr_\mathrm{\!ult}\{A_n\}.$
The class of subdirectly irreducible algebras in $\Pv$ that are
embeddable in members of a given $\Pr_\mathrm{\!ult}\{A_i\}$
is contained in the $\!1\!$-generator quasivariety
$\Su\Pr\Pr_\mathrm{\!ult}\{A_i\},$ hence by the implication
(\ref{d.qv1gen})$\Rightarrow$(\ref{d.setcompat}), each of these $n$
classes has the property that all its subsets are $\!\Pv\!$-compatible.
This gives~(\ref{d.*kdecomp}).

On the other hand, given any infinite $\kappa,$ let $\V$ be the
variety of sets given with a $\!\kappa\!$-tuple of zeroary operations
(constants) $c_\alpha$ $(\alpha\in\kappa).$
Since as a quasivariety, $\V$ is generated by its finitely
presented objects \cite[Proposition~2.1.18]{VAG},
and there are only $\kappa$ of these, it satisfies~(\ref{d.qv*k}).
On the other hand, since all operations are zeroary, every
equivalence relation on a $\!\V\!$-algebra is a congruence,
so the subdirectly irreducible algebras are precisely
the $\!2\!$-element algebras.
We have one of these for each $\!2\!$-class equivalence
relation on $\kappa,$ and one more corresponding to the
partition of that set into $\kappa$ and $\emptyset.$
This gives $2^{\kappa}$ subdirectly irreducible algebras,
no two of which are $\!\V\!$-compatible.
\end{proof}

\section{Afterthoughts on $\!\Pv\!$-compatible algebras.}\label{S.comfort}

Perhaps the concept of ``$\Pv\!$-compatible algebras''
is not the best handle on the phenomena we have been examining;
or at least should be complemented by another way of looking at them.
Suppose that for algebras $A$ and $B$ in $\Pv,$ we say that $A$
is ``comfortable'' with $B$ in $\Pv$ if the coprojection
$A\to A\cP_\Pv\,B$ is one-to-one; equivalently, if $A$ is
$\!\Pv\!$-compatible with some homomorphic image of $B$ in $\Pv.$
This relation is not in general symmetric; e.g., in the
variety of associative unital rings, each ring $\mathbb{Z}/n\mathbb{Z}$
is comfortable with $\mathbb{Z},$ but not vice versa.
Algebras $A$ and $B$ are $\!\Pv\!$-compatible if and only if each is
comfortable with the other.
(So the relations of $\!\Pv\!$-compatibility and of being
comfortable in $\Pv$ may be characterized in terms of one another.)
More generally, an arbitrary family of
algebras is $\!\Pv\!$-compatible if and only if each
is comfortable with the coproduct of the others.

If $\Pv$ is generated as a prevariety by a class of algebras $\X,$ we
see that an algebra $A$ is comfortable in $\Pv$ with an
algebra $B$ if and only if homomorphisms into members
of $\X$ that contain homomorphic images of $B$ separate points of $A.$
Hence, if we classify algebras $B\in\Pv$ according to which
algebras $A$ are comfortable with them,
then algebras $B_1$ and $B_2$ will belong to the same equivalence
class under this relation if the subclass of $\X$
consisting of algebras containing homomorphic images of $B_1$
coincides with the subclass of those containing images of $B_2.$
(We do not assert the converse.)
In particular, if $\Pv$ is residually small, so that
$\X$ can be taken to be a set, the number of these equivalence classes
has the cardinality of a set.
More generally, if $\Pv$ is generated by the union of $\kappa$ classes
of algebras, each absolutely directed under the relation $\preceq$ of
Theorem~\ref{T.absdird}, the same reasoning shows that it will have
at most $2^\kappa$ equivalence classes under this equivalence relation.
On the other hand, if we classify algebras according to which
other algebras they {\em are comfortable with}, we may, so far as I can
see, get up to $2^{2^\kappa}$ classes.

For any algebra $A$ in $\Pv,$ the class of algebras
which are comfortable with $A$ forms a subprevariety of $\Pv.$
The class of algebras that $A$ is comfortable with likewise yields
a subprevariety on throwing in the trivial algebra.
(A stronger statement, also easy to
see, is that this class is closed under
taking subalgebras and under taking products with arbitrary algebras
in $\Pv;$ equivalently, that if this class contains an algebra $B,$ then
it contains every algebra in $\Pv$ admitting a homomorphism to $B.)$

\section{On infinite symmetric groups: an answer and a question.}\label{S.Sym}
This last section does not depend on any of the preceding material.

It was shown in \cite{debruijn} (cf.~\cite{embed})
that for $\Omega$ an infinite set,
the group $S=\mathrm{Sym}(\Omega)$ of all permutations of $\Omega$
contains a coproduct of two copies of itself (from which it was
deduced by other properties of that group that it
contains a coproduct of $2^{\mathrm{card}(\Omega)}$ copies of itself).
In \cite[Question~4.4]{embed}, I asked, inter alia, whether, for every
subgroup $B$ of $S,$ if we regard $S$ as
a member of the variety of groups given with homomorphisms of $B$ into
them, $S$ contains a coproduct of two copies of itself
{\em in that variety}.

The answer is negative.
To see this, pick any $x\in\Omega$ and let $B$ be the
stabilizer in $S$ of $x.$
Writing elements of $S$ to the left of
their arguments and composing them accordingly, we see that
the partition of $S$ into left cosets of $B$
classifies elements according to where they send $x,$ and that
for each $y\in\Omega,$ the coset sending $x$ to $y$
has elements of finite order; e.g.,
if $y\neq x,$ the $\!2\!$-cycle interchanging $x$ and $y.$

On the other hand, I claim that if
$S_1$ and $S_2$ are any two groups with a common subgroup $B$
proper in each, then in the coproduct with amalgamation $S_1\cP_B\,S_2$
there are left cosets of $B$ containing no elements of finite order.
Indeed, the standard normal form in that coproduct shows that
each left coset is generated by a possibly empty alternating
string of left coset representatives of $B$ in $S_1$ and $S_2.$
When that string is nonempty and has even length,
one sees that elements of finite order cannot occur.
Hence for $S=\mathrm{Sym}(\Omega)$ and $B$ as above, $S,$
as a group containing $B,$ cannot contain a copy of $S_1\cP_B\,S_2.$

So let us modify our earlier question.

\begin{question}\label{Q.B<S}
For $\Omega$ an infinite set, what nice conditions, if any, on
a subgroup $B\subseteq S=\mathrm{Sym}(\Omega)$ will imply that
$S$ has a subgroup containing $B$ and
isomorphic over $B$ to $S\cP_B\,S$?

For instance, will this hold if $B$ is equal to,
or contained in, the stabilizer of a subset
of $\Omega$ having the same cardinality as $\Omega$?
If $B$ is finite?
\end{question}

In that same question in \cite{embed}, I asked whether
for any submonoid $B$ of the monoid $\mathrm{Self}(\Omega)$
of self-maps of $\Omega,$ the monoid $\mathrm{Self}(\Omega)$
must contain, over $B,$ a coproduct of two copies of itself with
amalgamation of $B.$
It seems likely that the subgroup
$B\subset\mathrm{Sym}(\Omega)\subset\mathrm{Self}(\Omega)$ used above
also gives a counterexample to this part of the question.
This will be so if we can show
that the subgroup of invertible elements of the monoid
coproduct of two copies of $\mathrm{Self}(\Omega)$ with amalgamation
of $B$ is isomorphic to the group coproduct of
two copies of $\mathrm{Sym}(\Omega)$ with amalgamation of $B,$
since we have seen that this is not embeddable over $B$ in the group
$\mathrm{Sym}(\Omega)$ of invertible
elements of $\mathrm{Self}(\Omega).$
But the analysis of coproducts of monoids with amalgamation,
even when the submonoid being amalgamated is a group, seems difficult.

The final part of that question posed the same problem for the
endomorphism algebra of an infinite-dimensional vector space
over a field.
To this I also do not know the answer.

\section{Glossary for the non-expert in universal algebra.}\label{S.glossary}

I indicate below the meanings of some basic concepts of
universal algebra, though more briefly and informally than would be
done in a textbook presentation.
(Definitions of some other concepts are recalled in the sections where
they are used.
I do not define concepts of category theory,
such as coproduct; or of set theory, such as ultraproduct,
or the distinction between sets and proper classes.
For these, see standard references such as
\cite{CW}, \cite{Ch+Keis}.)\vspace{.2em}

An {\em $\!n\!$-ary operation} on a set $X$ means a function
$X^n\to X;$ here $n$ is called the {\em arity} of the operation.
An {\em algebra} is a set given with a family of operations
of specified arities.
The list of operation-symbols and their
arities is the {\em type} of the algebra (used here
only in the phrase ``algebras of the same type'').
Constants in the definition of an algebra
structure (e.g., the $0$ and $1$ of a ring structure) are in this
note regarded as {\em zeroary} operations; indeed, $X^0$ is a
$\!1\!$-element set, so a map $X^0\to X$ specifies an element of $X.$
Given a subset $S$ of an algebra $X,$ the {\em subalgebra}
of $X$ generated by $S$ is here denoted $\langle S\rangle.$

The given operations of an algebra are called {\em primitive}
operations.
Expressions in a family of variable-symbols and iterated applications
of the primitive operations determine {\em derived operations}.
Such expressions are themselves called {\em terms}.
For instance, $(x\,y)\,z$ and $x\,(y\,z)$ are distinct
ternary terms in the operations of a group.
(They must be
considered distinct so that they can be used to write the
group identity of associativity.)
The variable-symbols are also considered terms;
they are the starting-point for the recursive construction
of all terms.
This technical sense of ``term'' will not stop us from
using the word in other ways,
e.g., in referring to the $\!m\!$-th term of a sequence.

An algebra all of whose primitive operations have finite
arity is called {\em finitary}.
(This does not preclude there being infinitely {\em many} primitive
operations; e.g., we have this for modules over an infinite ring.)

As indicated in \S\ref{S.prologue}, a {\em variety} of algebras
is the class of all algebras of a given type satisfying a
given set of identities.
In any variety $\V,$ one can construct a {\em free} algebra on any set,
satisfying the usual universal property.

The above concepts are assumed from \S\ref{S.free} on.
Starting with \S\ref{S.quasi}, we also refer to the
variety of algebras {\em generated} by
a family $\X$ of algebras of a given type, i.e., the
least variety containing $\X.$
This is clearly the class of all algebras that satisfy all
identities satisfied by all members of $\X.$
{\em Birkhoff's Theorem} states that it is also the class of
all homomorphic images of subalgebras of (generally infinite)
direct products of members of $\X,$ abbreviated
$\mathbb{H}\,\Su\Pr(\X).$
(Definitions of {\em prevariety} and {\em quasivariety},
and results for these analogous to Birkhoff's
Theorem, are recalled in~\S\ref{S.quasi}.)

To motivate a concept used from \S\ref{S.compat1} on, note that
if an algebra $A$ is {\em embedded} in a direct
product $\prod_I A_i,$ by a homomorphism with components
$f_i: A\to A_i,$ then $A\cong f(A)\subseteq\prod_I f_i(A).$
Modeled on the properties of this subalgebra,
one defines a {\em subdirect
product} of a family of algebras $(B_i)_{i\in I}$ to be a subalgebra
of $\prod_I B_i$ which projects surjectively to each $B_i.$
An algebra that, up to isomorphism,
cannot be so expressed without one of the projection
maps being an isomorphism is called {\em subdirectly irreducible}.


\end{document}